\documentclass[mnsc,nonblindrev]{informs3}

\OneAndAHalfSpacedXI

\usepackage[plainpages=false,hyperfootnotes=false]{hyperref}
\hypersetup{
  colorlinks   = true, 
  urlcolor     = blue, 
  linkcolor    = red, 
  citecolor   = blue 
}

\usepackage{natbib}
 \bibpunct[, ]{(}{)}{,}{a}{}{,}%
 \def\bibfont{\small}%
\TheoremsNumberedThrough     
\ECRepeatTheorems

\EquationsNumberedThrough

\MANUSCRIPTNO{}

\usepackage{dsfont}
\usepackage{amsmath}
\usepackage{amssymb}
\usepackage{enumitem}
\usepackage{booktabs}
\usepackage{siunitx}
\usepackage[colorinlistoftodos]{todonotes}
\usepackage{cleveref}
\usepackage{caption}
\usepackage{anyfontsize}
\usepackage{yhmath}
\usepackage{bm}
\usepackage{siunitx}

\usepackage{tabularx}
\newcolumntype{L}{>{\raggedright\arraybackslash}X}

\usepackage{tikz}
\usetikzlibrary{arrows.meta, decorations.pathmorphing}

\usepackage{algorithm}
\usepackage{algpseudocode}

\DeclareFontShape{OT1}{cmr}{bx}{sc}{<-> cmbcsc10}{}

\newcommand{\abs}[1]{\left\lvert #1 \right\rvert}
\newcommand{\norm}[1]{\left\lVert #1 \right\rVert}  

\newcommand{\I}{\mathcal{I}}

\newcommand{\B}{\mathcal{B}}

\newcommand{\K}{\mathcal{K}}

\DeclareRobustCommand{\L}{\ensuremath{\mathcal{L}}}

\newcommand{\ep}{{\varepsilon}}
\renewcommand{\epsilon}{{\varepsilon}}
\newcommand{\Id}{{\mathds{1}}}

\newcommand{\E}{{\mathbb{E}}}
\newcommand{\N}{{\mathbb{N}}}
\DeclareRobustCommand{\P}{\ensuremath{\mathbb{P}}}

\newcommand{\R}{{\mathds{R}}}

\newcommand{\x}{{\bm{x}}}
\newcommand{\z}{{\bm{z}}}

\newcommand{\X}{{\bm{X}}}

\newcommand{\Y}{{\bm{Y}}}
\newcommand{\y}{{\bm{y}}}
\newcommand{\Z}{{\bm{Z}}}

\newcommand{\bbeta}{{\bm{\beta}}}

\newcommand{\mN}{{\mathcal{N}}}

\DeclareMathOperator{\sgn}{sgn}

\newcommand{\eqd}{{\overset{d}{=}}}

\newcommand{\DB}{{\mathcal{B}}}

\newcommand{\iso}[1]{#1^{\uparrow}}
\newcommand{\Av}{\text{Av}}

\newcommand{\U}{{\mathcal{U}}}
\newcommand{\bU}{\ensuremath{{\bm{\mathcal{U}}}}}
\newcommand{\mfB}{\mathfrak{B}}
\newcommand{\bmfB}{\bm{\mathfrak{B}}}
\renewcommand{\S}{{\mathcal{S}}}
\newcommand{\tS}{{\tilde{\mathcal{S}}}}

\newcommand{\bphi}{\Phi}

\newcommand{\boldeta}{\bm{\eta}}
\newcommand{\bzeta}{\bm{\zeta}}

\newcommand{\overbar}[1]{\mkern 1.5mu\overline{\mkern-1.5mu#1\mkern-1.5mu}\mkern 1.5mu}

\renewcommand{\tilde}{\widetilde}
\renewcommand{\hat}{\widehat}

\usepackage[plainpages=false,hyperfootnotes=false]{hyperref}
\hypersetup{
  colorlinks   = true, 
  urlcolor     = blue, 
  linkcolor    = red, 
  citecolor   = blue 
}

\newenvironment{proof}[1][Proof]{\noindent\textbf{#1.} }{\hfill$\square$\par\medskip}

\begin{document}

\RUNAUTHOR{Tam and Pesenti}

\RUNTITLE{Bounds for DRO Problems}

\TITLE{Bounds for Distributionally Robust Optimization Problems}

\ARTICLEAUTHORS{%
\AUTHOR{Brandon Tam}
\AFF{Department of Statistical Sciences, University of Toronto \EMAIL{brandontam.tam@mail.utoronto.ca}} 
\AUTHOR{Silvana M. Pesenti}
\AFF{Department of Statistical Sciences, University of Toronto \EMAIL{silvana.pesenti@utoronto.ca}} 
\today\footnote{First version: April 8, 2025, with title ``Dimension reduction of distributionally robust optimization problems''.}
} 

\ABSTRACT{We study distributionally robust optimization (DRO) problems with uncertainty sets consisting of high-dimensional random vectors that are close in the multivariate Wasserstein distance to a reference random vector. We give conditions when the images of these sets under scalar-valued aggregation functions are contained in and contain uncertainty sets of univariate random variables defined via a univariate Wasserstein distance. This provides lower and upper bounds for the solution to general multivariate DRO problems that are computationally tractable. Furthermore, we generalize the results to uncertainty sets characterized by Bregman-Wasserstein divergences, which allows for asymmetric deviations from the reference random vector. Moreover, for DRO problems with risk measure criterion in the class of signed Choquet integrals, we derive semi-analytic formulae for the upper and lower bounds and the distribution that attains these bounds.}

\KEYWORDS{distributionally robust optimization, Wasserstein distance, Bregman-Wasserstein divergence, signed Choquet integral} 

\maketitle






\section{Introduction}\label{sec:intro}

Optimization problems with stochastic components (stochastic programs) occur frequently in the finance and risk management literature. The classical approach to handle stochasticity, proposed independently by \citet{beale1955minimizing} and \citet{dantzig}, is to optimize the expected value of the stochastic objective. This results in problems of the form
\begin{equation} \label{eq:stochastic_program}
\inf_{\boldsymbol{a}\in\mathcal{A}}\E[\ell(\boldsymbol{a},\X)],
\end{equation}
where $\mathcal{A}$ is a set of feasible actions, $\ell$ is a real-valued loss function, $\X$ is a random vector, and the expectation is taken with respect to (wrt) $\X$. Although $\ell$ is assumed to be affine in the original problem, subsequent developments in optimization theory; see, e.g., \citet{rockafellar1970convex, barbu2012convexity}, have allowed for the study of more general real-valued loss functions. We refer to \citet{shapiro2021lectures} for a comprehensive discussion of stochastic programs.

To evaluate the expectation in \eqref{eq:stochastic_program}, the distribution of $\X$ is assumed to be known. However, in many practical applications, the distribution is only partially known or needs to be estimated. Difficulties arising from estimation or partial knowledge of the underlying distribution are well studied in finance (\citet{michaud}), decision theory (\citet{smith}), and risk management (\citet{cont, embrechts, pesenti4}). To address this restriction, \citet{scarf1958min} proposed a distributionally robust version of \Cref{eq:stochastic_program}, termed distributionally robust optimization (DRO). DRO problems take the form 
\begin{equation} \label{eq:DRO}
\inf_{\boldsymbol{a}\in\mathcal{A}}\, 
    \sup_{\Z\in\mathcal{U}} \, \E[\ell(\boldsymbol{a},\Z)]\,,
\end{equation}
where $\mathcal{U}$ is a set of plausible alternative random vectors centered around the reference random vector $\X$--- so--called uncertainty sets. Over the past two decades, DRO problems have become popular in operations research and economics; indicatively see \citet{rahimian} and \citet{hansen}. 

Early research in DRO focuses on linear programs with ellipsoidal uncertainty sets (\citet{ben-tal, ben-tal2}). Here we work, among others, with uncertainty sets defined by a Wasserstein distance constraint. The Wasserstein distance is a well-studied metric originating from the field of optimal transport (\citet{villani}). This metric is now widely used outside of the theory of optimal transport, with important applications in mathematics, probability, and statistical theory (\citet{panaretos, dobrushin, munk}). In particular, worst-case expectation problems (i.e., the inner optimization problem in \eqref{eq:DRO}) with Wasserstein uncertainty are well studied in the operations research literature. For instance, \citet{zhao2018data} provide a characterization for the worst-case distribution (i.e., the distribution of the random vector attaining the inner supremum in \eqref{eq:DRO}) in data-driven risk-averse two-stage stochastic optimization problems when the reference distribution is the empirical distribution of observed data. In \citet{esfahani2018data}, the authors consider loss functions of the form $\ell(\x)=\max_{k\in\{1, \ldots K\}}\ell_k(\x)$, reformulate the worst-case expectation problem as a finite convex program, and obtain an explicit upper bound for the worst-case expectation when $\ell$ is proper, convex, and lower semi-continuous. For arbitrary loss functionals $\ell$, \citet{blanchet} and \citet{gao} derive a strong duality reformulation for the worst-case expectation problem. Another related Wasserstein DRO problem is the distributionally robust chance-constrained (DRCC) problem (\citet{xie2021distributionally, chen2024data}). DRCC problems take the form 
\begin{equation*}
\min_{\boldsymbol{a}\in\mathcal{A}}
    \left\{\boldsymbol{c}^T\boldsymbol{a}\right\} 
\ \text{subject to} \ 
    \inf_{\Z\in\mathcal{U}}\P\big(\Z\in\mathcal{C}({\boldsymbol{a}})\big)\geq 1-\epsilon,
\end{equation*}
where $\boldsymbol{a}, \boldsymbol{c}\in\R^m$ and $\mathcal{C}(\boldsymbol{a})$ is a constraint set that depends on the decision variable $\boldsymbol{a}$. In \citet{chen2024data}, the authors study the chance constraint and derive an exact reformulation for data-driven chance constraints. In contrast, \citet{xie2021distributionally} derives inner and outer approximations for the set of actions satisfying the chance constraint. Both \citet{xie2021distributionally} and \citet{chen2024data} assume that there exists $N$ historical samples from the true underlying distribution and the distribution of the reference random vector is the empirical distribution of these samples. In this paper, the reference distribution is not necessarily the empirical one.

Additionally, in our work, we consider both the Wasserstein distance and an asymmetric generalization of the Wasserstein distance known as the Bregman-Wasserstein (BW) divergence. The BW divergence was first introduced by \citet{carlier}, its geometry is discussed in \citet{kainth2025bregman}, and BW uncertainty sets are studied in \citet{guo,pesenti,pesenti3}, and \citet{pesenti5}. Here, we focus on the inner optimization problem in \eqref{eq:DRO} and assume that the risk factor $\X$ is multivariate. Specifically, we consider uncertainty sets defined using multivariate Wasserstein distance and multivariate BW divergence constraints. DRO problems involving the multivariate Wasserstein distance are challenging because explicit formulas for the multivariate Wasserstein distance between distributions are only known for a few distributions such as the Gaussian (\citet{dowson, gelbrich}). Numerical algorithms, e.g., the Sinkhorn algorithm, can be used in multivariate settings, but it is well known that the computational cost significantly increases with the dimension of $\X$ (\citet{cuturi}). We contribute to the literature by proposing upper and lower bounds to the inner optimization problem in \eqref{eq:DRO} that depend only on a univariate Wasserstein distance, regardless of the dimension of $\X$. We further derive expressions for the bounds that allow for efficient numerical evaluation, and characterize special cases under which a closed-form expression for the bounds exists. We also identify cases when the upper and lower bounds coincide and derive bounds for uncertainty sets characterized via the BW divergence. For some choices of BW uncertainty sets, we establish an intimate connection between the uncertainty set and the aggregation function via composable Bregman divergences, see \Cref{sec:composable}. 

We illustrate our lower and upper bounds on risk-aware DRO problems, where the expected value in \eqref{eq:DRO} is replaced by a law-invariant risk functional $\rho$ (also called a risk measure). Risk functionals are well studied in the literature. One of the earliest and most significant classes of risk functionals is the class of coherent risk measures, which was first introduced in \citet{artzner}. This class includes the expected value and the popular Expected Shortfall (ES). Over the past two decades, various generalizations have been proposed, including convex risk measures (\citet{follmer}), generalized deviation measures (\citet{rockafellar}), and the characterization of law-invariant coherent risk measures (\citet{kusuoka}). More recently, \citet{wang} provided an extensive study of signed Choquet integrals, which are a large class of risk functionals subsuming the well-known class of distortion risk measures. DRO problems with distortion risk measures are for example studied in \citet{cai}, \citet{Pesenti2023SIM}, \citet{bernard}, \citet{Coache2026WP}, and \citet{moresco2024uncertainty}. \citet{cai} consider moment based constraints for multivariate risks, \citet{Pesenti2023SIM} study portfolio optimization, \citet{bernard} focus on Wasserstein distance constraints for univariate risks, and \citet{Coache2026WP} and \citet{moresco2024uncertainty} work with robust dynamic risk measures. Another related work is \citet{bernard2}, who consider BW uncertainty around the choice of distortion risk measure rather than the underlying distribution. In our work, we consider multivariate risk factors, any law-invariant risk functional, and provide an application of our results to the class of signed Choquet integrals.

The paper is structured as follows. \Cref{sec:notation} introduces the notation, while \Cref{sec:Wasserstein} is devoted to the inner optimization problem in \eqref{eq:DRO} for law-invariant risk functionals when $\mathcal{U}$ is a multivariate Wasserstein uncertainty set and $\ell$ is a Lipschitz continuous function. In \Cref{sec:bw}, we generalize our results to the case when $\mathcal{U}$ is a BW uncertainty set. \Cref{sec:application} contains explicit bounds when the risk functional is a signed Choquet integral, and \Cref{sec:examples} illustrates numerical examples. 

\section{Multivariate Wasserstein Uncertainty}\label{sec:notation}

We work throughout with a non-atomic probability space $(\Omega, \B, \P)$, where $\B:=\B(\R^n)$ is the Borel sigma algebra on $\R^n$ for some $n\in\N^+$. We further denote vector-valued random variables (rvs) on $(\Omega, \B, \P)$ with capital boldface letters and real-valued vectors with lowercase boldface letters. For a random vector $\X=(X_1, \ldots, X_n)$, we write $F_{\X}(\x):=\P(\X\leq\x)$, $\x=(x_1, \ldots, x_n)\in\R^n$ for the cumulative distribution function (cdf) of $\X$, where vector inequalities are understood component-wise. For $a\geq 1$, we let $\norm{\x}_a:=(\sum_{i\in\mN}\abs{x_i}^a)^{\frac{1}{a}}$ denote the $\L^a$ norm on $\R^n$, where $\mN:=\{1,\ldots, n\}$. Furthermore, for $p\ge 1$ we denote by $\L^p_n$ the space of all $n$-dimensional random vectors $\X$ on $(\Omega, \B, \P)$ such that $\E[\norm{\X}_a^p]<\infty$, and by $\mathcal{M}_p(\R^n)$ the corresponding set of cdfs. As $a$ is fixed throughout, we write $\L^p_n$ and omit its dependence on the norm $a$. Moreover, if $n = 1$, we omit the dimension subscript and simply write $\L^p$. For any pair of random vectors $\X,\Y$, we denote by $\X\eqd \Y$ equality in distribution and by $\X=\Y$ $\P$-almost sure equality. Additionally, we denote by $\R^n_{\geq}$ the set of vectors in $\R^n$ with non-negative components.

The ($p$)-Wasserstein metric with norm $a$ on $\R^n$ (which we just call the Wasserstein distance) between $F\in\mathcal{M}_p(\R^n)$ and $G\in\mathcal{M}_p(\R^n)$ is defined as
\begin{equation*}
W^{n}(F,G):=\inf_{F_{\X}=F, F_{\Y}=G}
    (\E[\norm{\X-\Y}_{a}^{p}])^{\frac{1}{p}},
\end{equation*}
where the infimum is taken over all cdfs on $\R^{n}\times\R^{n}$ with $n$-dimensional marginal distributions $F$ and $G$. As we fix $a$ and $ p$ throughout the exposition, we simply write $W^n$, where $n$ indicates the dimension of the marginals, and omit the dependence on the power $p$ and the norm $a$.
For $n=1$, we omit the dimension superscript, i.e., we write $W(\cdot, \cdot):= W^1(\cdot, \cdot)$, and note that the norm $a$ becomes irrelevant. For $n = 1$, by the well-known fact that the infimum is attained by the comonotonic coupling (\citet{dall}), the 1-dimensional Wasserstein distance has representation
\begin{equation} \label{eq:quantile_rep_wass}
W(F,G):=\Big(\int_{0}^{1}
    \abs{F^{-1}(s)-G^{-1}(s)}^p\,ds\Big)^{\frac{1}{p}}\,,
\end{equation}
where $F^{-1}(\alpha):=\inf\{x\in\R~|~F(x)\geq\alpha\}$, $\alpha\in (0,1)$, is the (left-continuous) quantile function of $F$.

Next, we introduce the notation for uncertainty sets characterized by the Wasserstein distance.

\begin{definition} [Wasserstein Uncertainty Sets] \label{definition:uncertainty_rv} For $\epsilon\geq 0$, the Wasserstein uncertainty set with radius $\epsilon$ around the reference random vector $\X\in\L^p_n$ is
\begin{equation}\label{eq:d-dim_uncertainty_rv}
\bU_{\epsilon}^{n}(\X):=\{\Z\in\L^p_n~|~W^{n}
    (F_{\Z}, F_{\X})\leq\epsilon\}.
\end{equation}
\end{definition}
Similar to the Wasserstein distance, for $n=1$ we omit the dimension superscript and simply write $\U_{\epsilon}(X):=\bU_{\epsilon}^{1}(X)$. The parameter $\epsilon$ (also called tolerance distance) represents the magnitude of uncertainty, that is, the larger $\ep$, the larger the uncertainty.

Throughout the paper, we focus on the inner problem of \eqref{eq:DRO} and use $g \colon \R^n\to\R$ to denote an aggregation (or prediction) function that maps input risk factors $\X$ to a univariate output decision random variable $g(\X)$. The function $g$ should be interpreted as $\ell(a, \cdot)$ for a fixed decision $a$, where $\ell$ is the loss function in \eqref{eq:DRO}. Here, we assume that the aggregation function is given, e.g., estimated via statistical methods, and that the uncertainty stems from the risk factors $\X$. In an investment setting, for example, $g(\X)$ might represent the payoff of a portfolio, where $\X$ is the vector of asset prices subject to uncertainty and $g$ depends on the type of investments. The objective of our DRO problem is a law-invariant risk functional $\rho:\L^p\to\R$, subsuming the expected value, where  $\rho$ is law invariant if $X\eqd Y$ implies $\rho(X)=\rho(Y)$. Specifically, for $n>1$, we consider
\begin{equation}\label{eq:DRO-g-formulation}
\sup_{\Y\in \bU^{n}_{\epsilon}(\X)}\rho\big(g(\Y)\big)    \,.
\end{equation}

For optimization problem \eqref{eq:DRO-g-formulation} to be non-trivial, we assume throughout that the aggregation function $g(\cdot)$ is non-constant.
\begin{assumption}
    The aggregation function $g\colon \R^n \to \R$ is non-constant. 
\end{assumption}

There are different streams of literature that aim at simplifying robust DRO problems. \citet{Pesenti2024MP} for example derives necessary and sufficient conditions on the uncertainty set and the risk measure to recast a non-convex risk-aware DRO problem as a convex one. Here, the authors mostly mostly work with univariate uncertainty sets. In contrast, the uncertainty set in our problem is defined via the multivariate Wasserstein distance, making it difficult to solve without resorting to numerical algorithms such as the Sinkhorn algorithm. Thus, instead of explicitly calculating $\eqref{eq:DRO-g-formulation}$, we derive upper and lower bounds that are characterized by tractable univariate DRO problems of the form 
\begin{equation}\label{eq:DRO-g-univariate}
\sup_{Z\in \U_{\epsilon'}(g(\X))}\rho(Z)    \,.
\end{equation}
This DRO problem is typically more tractable as the uncertainty is around the univariate rv $g(\X)$, see for instance \citet{liu2022inf} for a closed-form expression when $\rho$ is a distortion risk measure. 

Comparing DRO problems \eqref{eq:DRO-g-formulation} and \eqref{eq:DRO-g-univariate}, we observe that the former accounts for uncertainty in the risk factors $\X$ (which then propagates to uncertainty in $g(\X)$), while the latter models uncertainty directly on the aggregate $g(\X)$. While optimization problem \eqref{eq:DRO-g-univariate} is over the univariate uncertainty set $\U_{\epsilon'}(g(\X))$, the corresponding aggregate univariate uncertainty set of DRO problem \eqref{eq:DRO-g-formulation} is $g\big(\bU_{\epsilon}^{n}(\X)\big):=\{g(\Z)~|~\Z\in\bU_{\epsilon}^{n}(\X)\}$, that is, each random vector in the uncertainty set $\Z \in \bU_{\epsilon}^{n}(\X)$ is mapped to a potential aggregate output $g(\Z)$. Therefore, \eqref{eq:DRO-g-formulation} is equivalent to 
\begin{equation*}
\sup_{\Y\in \bU^{n}_{\epsilon}(\X)}\rho\big(g(\Y)\big)  
=
\sup_{Z \in g(\bU_{\epsilon}^{n}(\X))}\rho(Z) \,.
\end{equation*}
In general the two univariate uncertainty sets are not the same, i.e.,  $g\big(\bU_{\epsilon}^{n}(\X)\big)\neq \U_{\epsilon}(g(\X))$. In \Cref{fig:set}, we illustrate how the uncertainty in DRO problem \eqref{eq:DRO-g-formulation} propagates from the multivariate set $\bU_{\epsilon}^{n}(\X)$ to the univariate aggregate set $g\big(\bU_{\epsilon}^{n}(\X)\big)$.

\smallskip

\begin{center}
\begin{tikzpicture}[line join=round,line cap=round]


\begin{scope}

\shade[
    inner color=red!10,
    outer color=red!35
]
(0,1) ellipse[x radius=2.8cm,y radius=1.35cm];

\draw[red!35,opacity=.45,line width=0.4pt]
(0,2.35)
.. controls (-0.15,1.6) and (-0.15,0.4) ..
(0,-0.35);

\draw[red!35,opacity=.45,line width=0.4pt]
(0,2.35)
.. controls (-1.5,1.9) and (-1.5,0.0) ..
(0,-0.35);

\draw[red!35,opacity=.45,line width=0.4pt]
(0,2.35)
.. controls (1.5,1.9) and (1.5,0.0) ..
(0,-0.35);

\draw[red!35,opacity=.45,line width=0.4pt]
(-2.65,0.8)
.. controls (-1.0,1.15) and (1.0,1.15) ..
(2.65,0.8);

\draw[red!35,opacity=.45,line width=0.4pt]
(-2.1,1.45)
.. controls (-0.8,1.7) and (0.8,1.7) ..
(2.1,1.45);

\draw[red!35,opacity=.45,line width=0.4pt]
(-2.1,0.25)
.. controls (-0.8,0.05) and (0.8,0.05) ..
(2.1,0.25);

\draw[red!35,opacity=.30,line width=0.3pt]
(-2.3,0.35)
.. controls (-0.5,0.75) and (1.0,1.1) ..
(2.0,0.9);

\node at (0,1) {$\Z$};

\node (e) at (0,1.75) {};

\node at (0.25,-1.15)
{$\bU^{n}_{\epsilon}(\X)=\{\Z\,|\,W^n(F_{\Z},F_{\X})\le\epsilon\}$};

\end{scope}


\begin{scope}[xshift=6cm]

\filldraw[
    fill=blue,
    fill opacity=.30
]
plot[smooth cycle,tension=.5] coordinates{
    (-0.25,0)
    (1,0.25)
    (2,0.75)
    (0.5,1.25)
    (0.25,2.75)
    (-1.75,1.75)
    (-2.25,0.75)
};

\node at (-0.75,1) {$g(\Z)$};

\node (f) at (-0.75,1.25) {};

\node at (0,-0.75)
{$g\big(\bU^{n}_{\epsilon}(\X)\big)$};

\end{scope}


\draw[
    <-,
    blue,
    thick,
    -{Stealth[length=4mm]}
]
(e.north) to[out=35,in=145] (f.north);

\node at (2.6,3.0) {$g$};

\end{tikzpicture}

\captionsetup{type=figure}
\caption{Visualization of how the multivariate uncertainty set $\bU^{n}_{\epsilon}(\X)$ (left) for the DRO problem \eqref{eq:DRO-g-formulation} cascades through to the univariate aggregate uncertainty set $g\big(\bU^{n}_{\epsilon}(\X)\big)$ (right).}
\label{fig:set}

\end{center}

\section{DRO with Wasserstein Uncertainty}\label{sec:Wasserstein}

In this section, we prove conditions on the aggregation function $g$ such that the image (under $g$) of the multivariate uncertainty set around $\X$ is contained in a univariate uncertainty set around $g(\X)$, i.e. $g\big(\bU_{\epsilon}^{n}(\X)\big)\subseteq \U_{\epsilon'}(g(\X))$. This provides upper bounds for the DRO problem \eqref{eq:DRO-g-formulation} that are univariate DRO problems of type \eqref{eq:DRO-g-univariate}. We furthermore characterize special cases where the set inclusion becomes a set equality, i.e., $g\big(\bU_{\epsilon}^{n}(\X)\big)=\U_{\epsilon'}(g(\X))$, and thus the DRO problems \eqref{eq:DRO-g-formulation} and \eqref{eq:DRO-g-univariate} coincide. Note that the tolerance distances for the two uncertainty sets are in general different.

\subsection{Lipschitz Aggregation Functions}\label{sec:Lipschitz}

We first consider the class of $K$-Lipschitz continuous functions. Recall that a function $f:\R^n\to\R$ is $K$-Lipschitz (wrt the $\L^a$ norm) if $\abs{f(\x)-f(\y)}\leq K\norm{\x-\y}_a$ for all $\x,\y\in\R^n$. Our first main result states that for a $K$-Lipschitz aggregation function $g$, the uncertainty set stemming from the risk factors is contained in the univariate uncertainty set around the aggregate $g(\X)$ with a tolerance distance of $K\epsilon$ -- the product of the Lipschitz constant and the tolerance distance $\ep$ of the multivariate uncertainty set. Lipschitz continuous functions are prevalent in portfolio optimization, important for fairness assessment (\citet{dwork2012fairness}), as well as ubiquitous in machine learning. Indeed, any fully connected network (FCN) and any convolutional neural network (CNN) with Lipschitz activation functions (e.g., ReLU, SoftPlus, sigmoid) is itself Lipschitz continuous. We refer to \citet{virmaux2018NIPS,Jordan2020NIPS}, and \citet{Kim2021proce} for details and discussions on Lipschitz constants of neural networks.

\begin{theorem} [Lipschitz aggregation]\label{theorem:Lipschitz} Suppose that $g:\R^n\to\R$ is $K$-Lipschitz. If $\X\in\L^p_n$ with $g(\X)\in\L^p$, then for any $\epsilon\geq 0$,
\begin{equation*}
g\big(\bU_{\epsilon}^{n}(\X)\big)\subseteq
    \U_{K\epsilon}(g(\X)).
\end{equation*}
\end{theorem}

\begin{proof}
Let $Z\in g\big(\bU_{\epsilon}^{n}(\X)\big)$. Then, by definition of $g\big(\bU_{\epsilon}^{n}(\X)\big)$, there exists a random vector $\Z$ such that $W^{n}(F_{\Z}, F_{\X})\leq\epsilon$ and $Z=g(\Z)$. Therefore, 
\begin{align*}
W(F_{g(\Z)}, F_{g(\X)}) &= \inf_{X'\eqd g(\X), \ Z'\eqd g(\Z)} 
    \E[\abs{X'-Z'}^p]^{\frac{1}{p}} \\
&= \inf_{g(\X')\eqd g(\X), \ g(\Z')\eqd g(\Z)} 
    \E\big[\abs{g(\X')-g(\Z')}^p\big]^{\frac{1}{p}} \\
&\leq \inf_{g(\X')\eqd g(\X), \ g(\Z')\eqd g(\Z)}  
    \E\big[K^p\norm{\X'-\Z'}_{a}^{p}\big]^{\frac{1}{p}} \\
&\leq \inf_{\X'\eqd \X, \ \Z'\eqd \Z} 
    \E\big[K^p\norm{\X'-\Z'}_{a}^{p}\big]^{\frac{1}{p}} \\
&=K\, W^{n}(F_{\Z}, F_{\X}) \\
&\leq K\,\epsilon.
\end{align*}
Hence, we conclude that $Z \in \U_{K\epsilon}(g(\X))$.
\end{proof}

The above theorem implies that if the aggregation function $g$ is Lipschitz, then the DRO problem \eqref{eq:DRO-g-formulation} can be upper bounded by a univariate DRO problem of type \eqref{eq:DRO-g-univariate} whose tolerance distance is scaled by the Lipschitz constant of $g$. 

Next, we discuss cases when the set inclusion in \Cref{theorem:Lipschitz} becomes a set equality, and thus the DRO problems \eqref{eq:DRO-g-formulation} and \eqref{eq:DRO-g-univariate} coincide. We first show that Lipschitz continuity alone is in general not sufficient. To see this, let $g(\x)=\sin(x_1)$ (which is Lipschitz with $K=1$), $\epsilon=3$, $\X\in\L^p_n$, and $Z=g(\X)+\epsilon$. Since $\E[\abs{Z-g(\X)}^p])^{\frac{1}{p}}=\epsilon$, it follows that $Z\in\U_{K\epsilon}(g(\X))$. Furthermore, for any random vector $\Y$, the support of $g(\Y)=\sin(Y_1)$ is a subset of $[-1,1]$. Since $Z=\sin(X_1)+3$ is supported on a subset of $[2,4]$, it follows that there does not exist a random vector $\Z$ such that $g(\Z)=Z$ and therefore $Z\notin g\big(\bU_{\epsilon}^{n}(\X)\big)$. 

To show set equality, we observe that any $K$-Lipschitz function $f\colon \R^n \to \R$ has decomposition
\begin{equation} \label{eq:decomposition}
f(\x)=\tilde{f}(\x^{(1)})+\bbeta^T\x^{(2)},
\end{equation}
where $\x=(\x^{(1)}, \x^{(2)})$ with $\x^{(1)}\in\R^m$, $\x^{(2)}\in\R^{n-m}$ for some $m\leq n$, $\bbeta\in\R^{n-m}$, and $\tilde{f}:\R^m\to\R$ is Lipschitz continuous. When $f$ has no linear components, then $f=\tilde{f}$ and $m=n$. In contrast, when $f$ has at least one linear component, then $\tilde{f}$ and $m$ are no longer unique as linear functions are also Lipschitz continuous. Throughout, we always choose the smallest possible $m$ to obtain the smallest Lipschitz constant for $\tilde{f}$.   

In an investment setting, when $\x$ is a vector of risky asset prices and $f$ is the portfolio payoff, $\bbeta^T\x^{(2)}$ represents the payoff from investing in $\beta_i$ units of asset $x^{(2)}_i$ for $i\in\{1, \ldots, n-m\}$. Moreover, $\tilde{f}(\x^{(1)})$ can be used to model non-linear payoffs that arise from investments in financial derivatives such as options.

Our next result characterizes a univariate uncertainty set around $g(\X)$ that is contained in the image of the multivariate uncertainty set $g\big(\bU_{\epsilon}^{n}(\X)\big)$, yielding a lower bound for optimization problem \eqref{eq:DRO-g-formulation}. The result states that the radius of the smaller set is proportional to the norm of the coefficient vector $\bbeta$ in the linear part of the decomposition \eqref{eq:decomposition}. Combining this result with \Cref{theorem:Lipschitz} allows us to identify when set equality holds. We refer to \citet{mao} for a similar result in the case when $g$ is linear.  

\begin{proposition} \label{prop:lin_comp} 
Let $g:\R^n\to\R$ be a Lipschitz continuous function (wrt the $\L^a$ norm) with decomposition \eqref{eq:decomposition}. If $\X\in\L^p_n$ with $g(\X)\in\L^p$, then for any $\epsilon\geq 0$,
\begin{equation*}
\U_{\norm{\bbeta}_b\epsilon}(g(\X))\subseteq 
    g\big(\bU_{\epsilon}^{n}(\X)\big),
\end{equation*}
where $\frac{1}{a}+\frac{1}{b}=1$. 
\end{proposition}

\begin{proof}
First, we note that if $\norm{\bbeta}_b=0$, then $\U_{\norm{\bbeta}_b\epsilon}(g(\X))=\{V~|~V\eqd g(\X)\}$, so the set containment is trivial. Thus, we assume wlog that $\norm{\bbeta}_b>0$. 

\smallskip

Let $Z\in\U_{\norm{\bbeta}_b\epsilon}(g(\X))$. We consider the cases of $a>1$ and $a=1$ separately. First, assume that $a>1$. Define the $(n-m)$-dimensional vector  
\begin{equation*}
\bbeta^{\frac{b}{a}} := 
    (\sgn(\beta_1)\abs{\beta_1}^{\frac{b}{a}}, 
\ldots, \sgn(\beta_{n-m})\abs{\beta_{n-m}}^{\frac{b}{a}}),
\end{equation*}
where $\sgn(\cdot)$ is the sign function. Let $\Z=(\Z^{(1)}, \Z^{(2)})\in\R^m\times\R^{n-m}$, where $\Z^{(1)}:=(X_1, \ldots, X_m)$ and 
\begin{equation*}
\Z^{(2)} := \Big(X_{m+1}+\tfrac{(\beta^{\frac{b}{a}})_1
    (Z-g(\X))}{\bbeta^T\bbeta^{\frac{b}{a}}}, \ldots, X_{n} + 
\tfrac{(\beta^{\frac{b}{a}})_{n-m}(Z-g(\X))}
    {\bbeta^T\bbeta^{\frac{b}{a}}}\Big),
\end{equation*} 
and $(\beta^{\frac{b}{a}})_i$ is the $i$-th component of $\bbeta^\frac{b}{a}$. It suffices to show that $g(\Z)=Z$ and $W^{n}(F_{\Z}, F_{\X})\leq\epsilon$. By definition of $\Z$, we have
\begin{align*}
g(\Z) &= \tilde{g}(X_1, \ldots, X_m) + 
    \bbeta^T\X^{(2)} + \bbeta^T\Big(
\tfrac{\bbeta^{\frac{b}{a}}(Z-g(\X))}
    {\bbeta^T\bbeta^{\frac{b}{a}}}\Big) \\
&= g(\X) + Z - g(\X) \\
&= Z.
\end{align*}
Furthermore,
\begin{align*}
W^{n}(F_{\Z}, F_{\X}) &=
\inf_{\X'\eqd\X, \Z'\eqd\Z}\big(
    \E\big[\norm{\Z'-\X'}_{a}^{p}\big]\big)^{\frac{1}{p}} \\
&\leq \big(\E[\norm{\Z-\X}_{a}^{p}]\big)^{\frac{1}{p}} \\
&= \Big(\E\big[\lVert\bbeta^{\frac{b}{a}}\rVert_{a}^{p}
    \tfrac{\abs{Z-g(\X)}^p}{\lvert\bbeta^T\bbeta^{\frac{b}{a}}
\rvert^p}\big]\Big)^{\frac{1}{p}} \\
&= \Big(\E\big[\norm{\bbeta}^{\frac{p\,b}{a}}_{b}
    \tfrac{\abs{Z-g(\X)}^p}{\norm{\bbeta}^{p\,b}_b}\big]
\Big)^{\frac{1}{p}} \\
&= \frac{\big(\E[\abs{Z-g(\X)}^p]
    \big)^{\frac{1}{p}}}{\norm{\bbeta}_b} \\
&\leq \frac{\norm{\bbeta}_b}{\norm{\bbeta}_b}\epsilon \\
&=\epsilon,
\end{align*}
where the third equality follows from the fact that $\lVert\bbeta^{\frac{b}{a}}\rVert_{a}^{p}=(\sum_{i=1}^{n-m}\abs{\beta_i}^b)^{\frac{p}{a}}=\norm{\bbeta}^{\frac{p\,b}{a}}_{b}$ and $\lvert\bbeta^T\bbeta^{\frac{b}{a}}\rvert^p=\abs{\sum_{i=1}^{n-m}\abs{\beta_i}^{1+\frac{b}{a}}}^p=\norm{\bbeta}_b^{p\,b}$, completing the proof for the case of $a>1$.

\smallskip

If $a=1$, then let $j^*\in\{1, \ldots, n-m\}$ such that $\abs{\beta_{j^*}}=\norm{\bbeta}_{\infty}$, where $\norm{\cdot}_\infty$ denotes the $\L^{\infty}$ norm, i.e., $\norm{\bbeta}_{\infty} := \max\{|\beta_1|, \ldots|\beta_{n-m}|\}$. Let $\Z$ be the random vector $\X$ with the $(m+j^*)$-th coordinate replaced by $X_{m+j^*}+\tfrac{\sgn(\beta_{j^*})(Z-g(\X))}{\norm{\bbeta}_{\infty}}$. Then, by similar arguments to the case of $a>1$, we have $g(\Z)=Z$ and $W^{n}(F_{\Z}, F_{\X})\leq\epsilon$. 
\end{proof}

Thus, for a $K$-Lipschitz function $g$ with decomposition \eqref{eq:decomposition},  \Cref{theorem:Lipschitz} and \Cref{prop:lin_comp} yield the following sequence of set inclusions of univariate uncertainty sets
\begin{equation}\label{eq:set-inclusion}
\U_{\norm{\bbeta}_b\epsilon}(g(\X)) \subseteq
    g\big(\bU_{\epsilon}^{n}(\X)\big) \subseteq 
\U_{K\epsilon}(g(\X))\,.    
\end{equation}
The smallest uncertainty set $\U_{\norm{\bbeta}_b\epsilon}$ has a tolerance distance proportional only to the linear term in the decomposition \eqref{eq:decomposition}, while the largest set $\U_{K\epsilon}(g(\X))$ has a tolerance distance scaled by the Lipschitz constant of $g$. Thus, if $K=\norm{\bbeta}_b$, then $g\big(\bU_{\epsilon}^{n}(\X)\big)=\U_{K\epsilon}(g(\X))$, and we obtain an equivalent reformulation for DRO problem \eqref{eq:DRO-g-formulation} in terms of  a DRO problem of type \eqref{eq:DRO-g-univariate}. For $K>\norm{\bbeta}_b$, we obtain better bounds when the Lipschitz constant $K$ is close to the norm of $\bbeta$ (i.e., $g$ is close to linear). We refer the reader to \Cref{sec:examples} for examples and a comparison of the bounds for different values of $K-\norm{\bbeta}_b$.

Note that the largest set is not affected by the choice of decomposition \eqref{eq:decomposition}. In contrast, the set $\U_{\norm{\bbeta}_b\epsilon}(g(\X))$ is the largest when we choose the largest possible value for $\norm{\bbeta}_b$ (or equivalently, the smallest possible value for $m$).

Next, we state two cases of set equality (which implies equality of DRO problems \eqref{eq:DRO-g-formulation} and \eqref{eq:DRO-g-univariate}) in \Cref{cor:set_equality}. For this, we first connect the Lipschitz constants of $f$ and $\tilde{f}$ in the decomposition \eqref{eq:decomposition}. 

\begin{lemma}\label{lemma:Lipschitz} Let $g:\R^n\to\R$ be a $K$-Lipschitz function (wrt the $\L^a$ norm) with decomposition \eqref{eq:decomposition}. Then, 
\begin{equation*}
K\leq\min\left\{L+\norm{\bbeta}_b, n^{\frac{1}{b}}
    \max\{L, \norm{\bbeta}_b\}\right\},
\end{equation*} 
where $L$ is the Lipschitz constant of $\tilde{g}$ (i.e., the non-linear part of $g$) and $\frac{1}{a}+\frac{1}{b}=1$. 
\end{lemma}

\begin{proof}
Let $g(\x)=\tilde{g}(\x^{(1)})+\bbeta^T\x^{(2)}$. Then,
\begin{align*}
\abs{g(\x)-g(\y)} &= \abs{\tilde{g}(\x^{(1)}) + 
    \bbeta^T\x^{(2)}-\tilde{g}(\y^{(1)})-\bbeta^T\y^{(2)}} \\
&\leq \abs{\tilde{g}(\x^{(1)}) - \tilde{g}(\y^{(1)})} + 
    \abs{\bbeta^T(\x^{(2)}-\y^{(2)})} \\
&\leq L\norm{\x^{(1)}-\y^{(1)}}_a + 
    \norm{\bbeta}_b\norm{\x^{(2)}-\y^{(2)}}_a, 
\end{align*}
where the second inequality follows from H\"{o}lder's inequality. To show that $K\leq L+\norm{\bbeta}_b$, observe that
\begin{align*}
L\norm{\x^{(1)}-\y^{(1)}}_a + 
    \norm{\bbeta}_b\norm{\x^{(2)} -\y^{(2)}}_a &\leq 
L\norm{\x-\y}_a + \norm{\bbeta}_b\norm{\x-\y}_a \\
&= (L+\norm{\bbeta}_b)\norm{\x-\y}_a. 
\end{align*}
Next, we show that $K\leq n^{\frac{1}{b}}\max\{L, \norm{\bbeta}_b\}$. From above we have
\begin{align*}
\abs{g(\x)-g(\y)} & \le L\norm{\x^{(1)}-\y^{(1)}}_a + 
    \norm{\bbeta}_b\norm{\x^{(2)}-\y^{(2)}}_a \\
&\leq L\norm{\x^{(1)} -\y^{(1)}}_1 + 
    \norm{\bbeta}_b\norm{\x^{(2)}-\y^{(2)}}_1 \\
&\leq \max\{L, \norm{\bbeta}_b\}(\norm{\x^{(1)} -
    \y^{(1)}}_1 + \norm{\x^{(2)}-\y^{(2)}}_1) \\
&= \max\{L, \norm{\bbeta}_b\}\norm{\x-\y}_1 \\
&\leq n^{\frac{1}{b}}\max\{L, \norm{\bbeta}_b\}\norm{\x-\y}_a,
\end{align*}
where the first inequality follows from the fact that $\norm{\cdot}_a\leq\norm{\cdot}_1$ for $a\geq 1$ and the last inequality follows from the fact that $\norm{\cdot}_1\leq n^{\frac{1}{b}}\norm{\cdot}_a$. 
\end{proof}

\begin{corollary} \label{cor:set_equality} 
Suppose $g:\R^n\to\R$ satisfies the conditions of \Cref{lemma:Lipschitz}. Let $\X\in\L^p_n$ with $g(\X)\in\L^p$ and $\epsilon\geq 0$. Assume any one of the following holds:
\begin{enumerate}[label = $\roman*)$]
    \item $a=1$ and $\norm{\bbeta}_{\infty}\geq L$
    \item $m=0$
\end{enumerate}  
Then, $g\big(\bU_{\epsilon}^{n}(\X)\big)=\U_{K\epsilon}(g(\X))$.
\end{corollary}

\begin{proof}
By \eqref{eq:set-inclusion} it suffices to show that $K\leq\norm{\bbeta}_b$ for both cases. 

\smallskip

For $i)$, by \Cref{lemma:Lipschitz}, $K\leq\min\left\{L+\norm{\bbeta}_b, n^{\frac{1}{b}}\max\{L, \norm{\bbeta}_b\}\right\}=\min\{L+\norm{\bbeta}_{\infty}, \norm{\bbeta}_{\infty}\}=\norm{\bbeta}_{\infty}$.

\smallskip

For $ii)$, by \Cref{lemma:Lipschitz}, $K\leq\min\left\{L+\norm{\bbeta}_b, n^{\frac{1}{b}}\max\{L, \norm{\bbeta}_b\}\right\}=\min\{\norm{\bbeta}_{b}, n^{\frac{1}{b}}\norm{\bbeta}_{b}\}=\norm{\bbeta}_{b}$.
\end{proof}

Note that case $ii)$ corresponds to the case when $g$ is linear, which was first proven in \citet{mao}. Our result shows that equality also holds for a larger class of aggregation functions when the norm of the multivariate Wasserstein distance used to define the uncertainty set is 1. With the set inclusions \eqref{eq:set-inclusion} at hand, we can apply them to bound the risk-aware DRO problem \eqref{eq:DRO-g-formulation}.

\begin{proposition} [Worst-case risks under Lipschitz aggregation]\label{prop:Lipschitz}
Let $g:\R^n\to\R$ be a $K$-Lipschitz function with decomposition \eqref{eq:decomposition}. Further, let $\X\in\L^p_n$ with $g(\X)\in\L^p$, $\epsilon\geq 0$, and $\rho$ be a law-invariant risk functional. Then,
\begin{equation} \label{eq:Lipschitz_bound} 
\sup_{Y\in \U_{\norm{\bbeta}_b\epsilon}(g(\X))} 
    \rho(Y) \leq \sup_{\Y\in \bU^{n}_{\epsilon}(\X)} 
\rho\big(g(\Y)\big)\leq
    \sup_{Y\in \U_{K\epsilon}(g(\X))}\rho(Y)\,.
\end{equation}
Moreover, if $a=1$ and $\norm{\bbeta}_{\infty}\geq L$, or if $a\geq 1$ and $m=0$, then the inequalities become equalities. 
\end{proposition}

\begin{proof}
The inequalities follow from \Cref{theorem:Lipschitz} and \Cref{prop:lin_comp} and the equality cases follow from \Cref{cor:set_equality}. 
\end{proof}

One might wonder if scaling the aggregation function $g$ to obtain a smaller Lipschitz constant has any impact on the bounds. By positive homogeneity of the Wasserstein distance, it holds that $W(F_{c\,Y}, F_{c\,Z})=c\,W(F_Y, F_Z)$ for any $c>0$ and $Y, Z\in\L^p$. This implies that $\U_{c\,\epsilon}(g(\X))=\{Y~|~\tfrac{1}{c}\,Y\in\U_{\epsilon}(\tfrac{g(\X)}{c})\}$. Therefore, if $\rho$ is positive homogeneous (i.e., $\rho(\theta\,Y)=\theta\,\rho(Y)$ for all $\theta\geq 0$), then the upper bound in \eqref{eq:Lipschitz_bound} is
\begin{equation*}
 \sup_{Y\in \U_{K\epsilon}(g(\X))}\rho(Y) \,
=
\sup_{\left\{Y~\big|~\tfrac{1}{K}\,Y\in\U_{\epsilon}\big(\tfrac{g(\X)}{K}\big)\right\}}\rho(K\,\tfrac{1}{K}Y) = K\sup_{Z\in\U_{\epsilon}\big(\tfrac{g(\X)}{K}\big)}\rho(Z).
\end{equation*}
That is, the upper bound in \eqref{eq:Lipschitz_bound} is equal to $K$ times the worst-case risk when the uncertainty set has radius $\epsilon$ around the reference rv $\hat{g}(\X)$, where $\hat{g}(\cdot):=\frac{g(\cdot)}{K}$ is a $1$-Lipschitz function. 

\subsection{Locally Lipschitz Aggregation Functions}\label{sec:loc_Lipschitz}

We extend the results of the preceding section to $g$ being locally Lipschitz. Specifically we allow the Lipschitz condition on $\tilde{g}$ to hold locally instead of globally. That is, we work with aggregation functions $g:\R^n\to\R$ satisfying 
\begin{equation}\label{eq:g-decomposition-local}
g(\x)=\tilde{g}(\x^{(1)})+\bbeta^T\x^{(2)}\,,    
\end{equation}
where $\x=(\x^{(1)}, \x^{(2)})$ with $\x^{(1)}\in\R^m$, $\x^{(2)}\in\R^{n-m}$ for some $m\leq n$, $\bbeta\in\R^{n-m}$, and $\tilde{g}:\R^m\to\R$ is locally Lipschitz. Recall that a function $f\colon \R^n \to \R$ is locally Lipschitz if for any $\z\in\R^n$, there exists a neighbourhood $S$ of $\z$ such that $f$ restricted to $S$ is Lipschitz. 

In the locally Lipschitz case, we require some additional boundedness assumptions on the risk factors $\X:=(\X^{(1)}, \X^{(2)})$ with $\X^{(1)} \in \L^p_m$ and $\X^{(2)} \in \L^p_{n-m}$. In particular, we require the existence of a compact subset $C\subseteq\R^{m}$, such that 
\begin{equation}\label{eq:boundedness}
\text{supp}(\X^{(1)})\subseteq C,
\end{equation} 
where $\text{supp}(\Y)$ denotes the support of a random vector $\Y$. Note that the risk factors associated with the linear part of $g$ can be unbounded as long as $\X\in\L^p_n$ and $g(\X)\in\L^p$. Due to this boundedness assumption, we slightly modify the notation of the uncertainty sets used in \Cref{sec:Lipschitz}. 

\begin{definition} [Wasserstein Uncertainty] \label{definition:uncertainty_c} For $\epsilon\geq 0$, $m\leq n$, and $C\subseteq\R^m$, the Wasserstein uncertainty set with radius $\epsilon$ for random vectors satisfying \eqref{eq:boundedness} is 
\begin{equation}\label{eq:set_loc_Lip}
\bU^{n,m}_{\epsilon,C}(\X) := 
    \{\Z\in\L^p_n~|~W^{n}(F_{\Z}, F_{\X})\leq\epsilon 
\ \text{and}\ 
    \text{supp}\big((Z_1, \ldots, Z_{m})\big)\subseteq C\}, 
\end{equation}
where $\X=(\X^{(1)}, \X^{(2)})\in\L_p^n$ and $\text{supp}(\X^{(1)})\subseteq C$.
\end{definition}

We note that the set $\bU_{\epsilon}^{n}(\X)$ introduced in \Cref{sec:notation}  is a special case of \eqref{eq:set_loc_Lip} with $C=\R^m$. 

Before stating the next result, we recall a useful property of locally Lipschitz functions from \citet{scanlon}. This property guarantees the existence of the Lipschitz constant $K_C$ in the statement of \Cref{theorem:loc_lip}. 

\begin{lemma} \label{lemma:loc_lip} (Theorem 2.1 from \citet{scanlon}) Let $(M,d)$ be a metric space. A function $f:M\to\R$ is locally Lipschitz if and only if it is Lipschitz on each compact subset of $M$. 
\end{lemma}

The following result is the analogue to \Cref{theorem:Lipschitz} and \Cref{prop:lin_comp} for locally Lipschitz aggregation functions.

\begin{theorem} [Locally Lipschitz aggregation]\label{theorem:loc_lip} Suppose that $g:\R^n\to\R$ is locally Lipschitz with decomposition \eqref{eq:g-decomposition-local}. If $\X\in\L^p_n$ with $g(\X)\in\L^p$ and $C\subseteq\R^m$ is a compact set such that \eqref{eq:boundedness} holds, then for any $\epsilon\geq 0$,
\begin{equation*}
\U_{\norm{\bbeta}_b\epsilon}(g(\X)) \subseteq 
    g\big(\bU_{\epsilon,C}^{n,m}(\X)\big) \subseteq 
\U_{K_C\epsilon}\big(g(\X)\big),
\end{equation*}
where $K_C$ is the Lipschitz constant of $g$ restricted to the compact subset $C\subseteq\R^m$ and $\frac{1}{a}+\frac{1}{b}=1$. 
Moreover, $K_C\leq\min\left\{L_C+\norm{\bbeta}_b, n^{\frac{1}{b}}\max\{L_C, \norm{\bbeta}_b\}\right\}$, where $L_C$ is the Lipschitz constant of $\tilde{g}$ restricted to $C\subseteq\R^m$. 
\end{theorem}

\begin{proof}
We first show that $g\big(\bU_{\epsilon,C}^{n,m}(\X)\big)\subseteq\U_{K_C\epsilon}\big(g(\X)\big)$. Let $Z\in g\big(\bU_{\epsilon,C}^{n,m}(\X)\big)$. By definition of $g\big(\bU_{\epsilon,C}^{n,m}(\X)\big)$, there exists a random vector $\Z$ such that $W^{n}(F_{\Z}, F_{\X})\leq\epsilon$, $Z=g(\Z)$, and $\text{supp}(Z_1, \ldots, Z_m)\subseteq C$. Therefore, 
\begin{align*}
W(F_{g(\Z)}, F_{g(\X)}) &= 
    \inf_{X'\eqd g(\X), \ Z'\eqd g(\Z)}
\E[\abs{X'-Z'}^p]^{\frac{1}{p}} \\
&\leq \inf_{\X'\eqd\X, \ \Z'\eqd\Z}
    \E\big[\abs{g(\X')-g(\Z')}^p\big]^{\frac{1}{p}} \\
&\leq \inf_{\X'\eqd\X, \ \Z'\eqd\Z} 
    \E\big[K_C^p\norm{\X'-\Z'}_{a}^{p}\big]^{\frac{1}{p}} \\
&\leq K_C\,\epsilon.
\end{align*}
Hence, we conclude that $Z \in \U_{K_C\epsilon}(g(\X))$.

\smallskip

Next, we show that $\U_{\norm{\bbeta}_b\epsilon}(g(\X))\subseteq g\big(\bU_{\epsilon,C}^{n,m}(\X)\big)$. Again, we assume wlog that $\norm{\bbeta}_b>0$. Let $Z\in\U_{\norm{\bbeta}_b\epsilon}(g(\X))$. Next, define $\Z$ as in the proof of \Cref{prop:lin_comp}. It suffices to show that $g(\Z)=Z$, $W^{n}(F_{\Z}, F_{\X})\leq\epsilon$, and $\text{supp}(Z_1, \ldots, Z_m)\subseteq C$. The first two claims follows from the arguments in the proof of \Cref{prop:lin_comp}. Furthermore, since $(Z_1, \ldots, Z_m)=(X_1, \ldots, X_m)$, the last claim holds.

\smallskip

Finally, the bound on the constant $K_C$ follows from the same arguments as in the proof of \Cref{lemma:Lipschitz}. 
\end{proof}

Next, we state an analogous result to \Cref{cor:set_equality} for locally Lipschitz $g$. Note that $m=0$ (the second case of \Cref{cor:set_equality}) implies that $g$ is linear and hence globally Lipschitz. Therefore, we only state the analog for  case $i)$ of \Cref{cor:set_equality}. 

\begin{corollary} \label{cor:set_equality_loc_lip} 
Suppose $g:\R^n\to\R$ is locally Lipschitz with decomposition \eqref{eq:g-decomposition-local}. Let $\epsilon\geq 0$, $\X\in\L^p_n$ with $g(\X)\in\L^p$, and $C\subseteq\R^m$ be a compact set such that \eqref{eq:boundedness} holds. If $a=1$ and $\norm{\bbeta}_{\infty}\geq L_C$, then $g\big(\bU_{\epsilon,C}^{n,m}(\X)\big)=\U_{K_C\epsilon}(g(\X))$, where $K_C$ and $L_C$ are as defined in \Cref{theorem:loc_lip}.
\end{corollary}

Applying the above results to risk-aware DRO problems, we obtain the following bounds for \eqref{eq:DRO-g-formulation}.
\begin{proposition} [DRO under locally Lipschitz aggregation]\label{prop:loc_lip} Assume that the conditions of \Cref{theorem:loc_lip} are satisfied. Then, for any $\epsilon\geq 0$ and law-invariant risk functional $\rho:\L^p\to\R$,
\begin{equation*}
\sup_{Y\in \U_{\norm{\bbeta}_b\epsilon}(g(\X))} 
    \rho(Y) \leq 
\sup_{\Y\in \bU_{\epsilon,C}^{n,m}(\X)}\rho\big(g(\Y)\big)
    \leq\sup_{Y\in \U_{K_C\epsilon}(g(\X))}\rho(Y)\,.
\end{equation*}
Moreover, if $a=1$ and $\norm{\bbeta}_b\geq L_C$, then the inequalities become equalities.
\end{proposition}

\begin{proof}
This follows from \Cref{theorem:loc_lip}.
\end{proof}

\section{DRO with Bregman-Wasserstein Uncertainty} \label{sec:bw}

Here, we consider uncertainty sets defined via a BW divergence, which is a generalization of the Wasserstein distance. Differently to the Wasserstein distance, the Bregman divergence is asymmetric and thus suitable when uncertainties arise in an asymmetric manner around the reference distribution. Specifically, we consider uncertainty sets defined as all random vectors that have a BW divergence of at most $\ep$ around a reference random vector. We refer to \cite{pesenti3} for a portfolio application. 

Key results are in \Cref{sec:separable}, where we generalize \Cref{theorem:Lipschitz} to BW divergences with separable Bregman generators. In \Cref{sec:m_distance}, we study a special case of the Mahalanobis distance, which is an example of a BW divergence with a separable generator. For the Mahalanobis distance, we derive bounds for the worst-case risk under weaker assumptions on the risk functional. Finally, in \Cref{sec:composable}, we consider Bregman generators that are compositions of a univariate generator and a real-valued aggregation function.

\subsection{Multivariate Bregman-Wasserstein Uncertainty} \label{sec:notation_bregman}

Before defining the BW divergence, we recall the definitions of a Bregman generator and the Bregman divergence. The Bregman generator associated with a BW divergence determines its key properties. In particular, unlike the Wasserstein distance, the BW divergence is not symmetric unless its Bregman generator is of the form $\bphi(\x)=\x^TA\x$ for some positive definite matrix $A$. This asymmetry is useful, for example, if $g(\X)$ is the portfolio payoff of an investment. In this case, deviations in the negative direction (i.e., losses) should have larger weight than deviations in the positive direction (i.e., gains). 

\begin{definition}[Bregman Divergence]  \label{definition:bregman}
Let $\mathcal{X}^n\subseteq\R^n$ be a convex domain. A Bregman generator is a function $\bphi: \mathcal{X}^n\to\R, \ n\geq 1$, that is convex and differentiable. The Bregman divergence associated with Bregman generator $\bphi$ is defined as
\begin{equation*}
B_{\bphi}\big(\z_1, \z_2\big) := 
    \bphi(\z_1) - \bphi(\z_2) - 
\nabla\bphi(\z_2)\cdot(\z_1-\z_2)\,,\quad 
    \z_1,\z_2\in\mathcal{X}^n\,,
\end{equation*}
where $\nabla\bphi(\z)$ denotes the gradient of $\bphi$ and $\cdot$ denotes the dot product. For $n=1$, we use $\phi$ for univariate generators. 
\end{definition}

Note that the convexity of $\bphi$ guarantees that the Bregman divergence is non-negative. Moreover, if $\bphi$ is strictly convex, then the Bregman divergence $B_{\bphi}\big(\z_1, \z_2\big)$ is zero if and only if $\z_1=\z_2$ (and therefore is a mathematical divergence). Here, we use convexity in the non-strict sense unless otherwise stated, and we refer to \citet{pesenti3} for examples of non-strictly convex generators to model distributional uncertainty. 

\begin{definition} [BW Divergence] \label{definition:bw} Let $\bphi : \mathcal{X}^n \to \R$ be a Bregman generator. Then, the BW divergence associated with $\bphi$, from cdf $F$ to cdf $G$, is defined as
\begin{equation*}
\DB^{n}_{\bphi} [F, G] := 
    \inf_{F_{\X}=F, F_{\Y}=G}
\E\big[B_{\bphi}\big(\X, \Y\big)\big].
\end{equation*}
\end{definition}
For random vectors $\X, \Y$, with cdfs $F, G$, respectively, we use the notation $\DB^{n}_{\bphi} (\X, \Y) := \DB^{n}_{\bphi} [F_{\X}, F_{\Y}]$, that is, round brackets for random vectors and square brackets for cdfs. Moreover, we omit the dimension superscripts for $n=1$ and simply write $\DB_{\phi} [F, G]:= \DB^{1}_{\bphi} [F, G]$ and $\DB_{\phi} (\X, \Y) := \DB^{1}_{\bphi} (\X, \Y)$. 

When the Bregman generator is the squared 2-norm, i.e. $\bphi(\x) =\norm{\x}_2^2$, then the BW divergence reduces to the squared $2$-Wasserstein distance with norm $\L^2$. Furthermore, for $n=1$, the infimum is achieved by the comonotonic coupling (\citet{pesenti}). In other words, we can rewrite the 1-dimensional BW divergence as 
\begin{equation} \label{eq:1d_bregman}
\DB_{\phi}[F, G] = 
    \int_{0}^{1}B_{\phi}\big(F^{-1}(t), G^{-1}(t)\big)\,dt.
\end{equation}
For $n= 1$ and using \eqref{eq:1d_bregman}, it is straightforward to prove that the BW divergence is convex in its first component on the space of quantile functions (\citet{pesenti3}). This is however, not necessarily the case for arbitrary dimensions. Indeed, even though the Bregman divergence is convex in its first argument, the BW divergence is generally not convex in its first argument on the space of random vectors. We refer to \citet{pesenti}, \citet{kainth2025bregman}, and the references therein, for detailed discussions of the BW divergence.

With the definition of the BW divergence at hand, we now introduce the BW uncertainty set.

\begin{definition}[BW Uncertainty Sets] \label{definition:bw_set}
For $\epsilon\geq 0$, we define the BW uncertainty set around the reference random vector $\X$ associated with the Bregman generator $\bphi:\mathcal{X}^n\to\R$ by 
\begin{equation*}
\bmfB^{n}_{\bphi, \epsilon}(\X) := 
    \{\Z~|~\DB^{n}_{\bphi}[F_{\Z},F_{\X}]\leq\epsilon\}.
\end{equation*}
\end{definition}
Again, we omit the dimension superscript for $n=1$ and write $\mfB_{\phi, \epsilon}(X):=\bmfB_{\bphi, \epsilon}^{1}(X)$. Moreover, whenever we write $\mfB_{\phi, \epsilon}(X)$ or $\bmfB^{n}_{\bphi, \epsilon}(\X)$, we tacitly assume that the uncertainty sets contain at least two, and thus infinitely many elements.

\subsection{Separable Bregman Generators}\label{sec:separable}

In this section, we consider separable Bregman generators, which are functions $\bphi:\mathcal{X}^n\to\R$ of the form $\bphi(\x)=\sum_{k\in\mN}\phi_k(x_k)$, where $\phi_k:\mathcal{X}\to\R$, $k\in\mN$, are themselves Bregman generators. Note that the popular Kullback-Leibler (KL) divergence and the Mahalanobis distance are Bregman divergences with separable generators. We refer to \citet{kainth2025bregman} for a detailed discussion of the KL divergence and \Cref{sec:m_distance} for a discussion of the Mahalanobis distance. 

\begin{example}[KL Divergence]
{\parindent0pt Define 
\begin{equation*}
\S := \left\{\y\in\R^n~|~
    \y = \Big(\frac{e^{x_{1}}}{1+\sum_{i\in\mN}e^{x_i}},
\ldots, \frac{e^{x_{n}}}{1+\sum_{i\in\mN}e^{x_i}}\Big), 
    \,  \x\in\R^n\right\}.
\end{equation*} 
Then, $\Delta^{n}:=\{\z=(z_0,\y)\in\R^{n+1}~|~\y\in \S, \ z_0=1-\sum_{i\in\mN}y_i\}$ is the open unit simplex in $\R^{n+1}$. Consider the Bregman generator $\bphi(\z)=\sum_{k=0}^{n}z_k\log(z_k), \ \z\in\Delta^n$ (i.e., the negative Shannon entropy). Then, $\nabla\bphi(\z)=(\log(z_0)+1, \ldots, \log(z_n)+1)$ and thus
\begin{equation*}
B_{\bphi}[\z_1,\z_2] = 
    \sum_{i=0}^{n}z_{1i}\log\big(\frac{z_{1i}}{z_{2i}}\big),
\end{equation*}
which is the KL divergence on the simplex. If we remove the restriction onto the simplex, then we recover the generalized KL divergence (see e.g., \citet{miller}).  
}
\end{example}

In the next theorem, we establish a connection between the BW divergences for rvs and random vectors when the multivariate generator $\bphi$ is separable. In particular, we construct a new uncertainty set defined via univariate BW divergences that contains the multivariate uncertainty set $\bmfB_{\bphi, \epsilon}^{n}(\X)$. This connection allows us to derive an analogue to \Cref{prop:Lipschitz} for separable Bregman generators.  

\begin{theorem} [Separable Bregman generators]\label{theorem:bregman_sum} Let $\bphi(\x)=\sum_{k\in\mN}\phi_k(x_k)$ be a separable Bregman generator. Then, for any $\epsilon\geq 0$ and random vector $\X$,
\begin{enumerate}[label = $\roman*)$]
    \item $\bmfB_{\bphi, \epsilon}^{n}(\X)\subseteq \{\Z~|~\sum_{i\in\mN}\DB_{\phi_i}(Z_i, X_i)\leq\epsilon\}$
    \item $\bmfB_{\bphi, \epsilon}^{n}(\X) \cap \aleph
    \, = \{\Z~|~\sum_{i\in\mN}\DB_{\phi_i}(Z_i, X_i)\leq\epsilon\} \cap \aleph$,
where 
\begin{equation*}
\aleph := \{\Z~|~ Z_i, X_i 
    \ \text{are comonotonic\footnotemark{} for all} \ i\in\mN\}.
\end{equation*}
\end{enumerate}
\footnotetext{Two rvs $X, Y$ on $(\Omega, \B, \P)$ are comonotonic if there exists $\Omega_0\subseteq\Omega$ with $\P(\Omega_0)=1$ such that for all $\omega, \omega'\in\Omega_0$, $\big(X(\omega)-X(\omega')\big)\big(Y(\omega)-Y(\omega')\big)\geq 0$.}
\end{theorem}

\begin{proof}
We first prove $i)$. For any random vectors $\X$ and $\Z$, define the sets $\tS_i := \{(Z_i',X_i')~|~X_i'\eqd X_i,\; Z_i'\eqd Z_i\}$, $i\in\mN$, and
\begin{equation*}
\tS := \{(\Z',\X')~|~(X_i', Z_i') \in \tS_i
    \ \text{for all} \ i\in\mN\}.
    \quad \text{and} \quad
\S := \{(\Z',\X')~|~\X'\eqd\X, \Z'\eqd\Z\} \end{equation*}
Since $\bphi$ is separable,
\begin{align*}
\DB_{\bphi}^{n}(\Z, \X) &= \inf_{(\Z', \X')\in\S} 
    \E\big[\sum_{i\in\mN}\phi_i(Z_i') - 
\phi_i(X_i')-\phi_i'(X_i')(Z_i'-X_i')\big] \\
&= \inf_{(\Z', \X')\in\S}\, \sum_{i\in\mN}
    \E\big[\phi_i(Z_i')-\phi_i(X_i') - 
\phi_i'(X_i')(Z_i'-X_i')\big] \\
&\geq \inf_{(\Z', \X')\in\tS}\, 
    \sum_{i\in\mN}\E\big[\phi_i(Z_i') - 
\phi_i(X_i')-\phi_i'(X_i')(Z_i'-X_i')\big] \\
&= \sum_{i\in\mN}\inf_{(Z_i', X_i')\in\tS_i} 
    \E\big[\phi_i(Z_i')-\phi_i(X_i') - 
\phi_i'(X_i')(Z_i'-X_i')\big] \\
&=\sum_{i\in\mN}\DB_{\phi_i}(Z_i, X_i),
\end{align*}
where the inequality holds since $\S\subseteq\tS$. Moreover, the second last equality holds since $(\Z', \X')\in\tS$ if and only if $(Z_i', X_i')\in\tS_i$ for all $i\in\mN$ and each summand in the optimization objective depends only on one $i\in\mN$. 

\smallskip

Since $\sum_{i\in\mN}\DB_{\phi_i}(Z_i, X_i)\leq \DB_{\bphi}^{n}(\Z, \X)$, it follows that for any $\epsilon\geq 0$, the inequality $\DB_{\bphi}^{n}(\Z, \X)\leq\epsilon$ implies that $\sum_{i\in\mN}\DB_{\phi_i}(Z_i, X_i)\leq\epsilon$. Therefore, we conclude that $i)$ holds. 

\smallskip

Second, we prove $ii)$. From $i)$, it follows that 
\begin{equation*}
\bmfB_{\bphi, \epsilon}^{n}(\X) \cap \aleph
\subseteq \{\Z~|~\sum_{i\in\mN}\DB_{\phi_i}
    (Z_i, X_i)\leq\epsilon\} \cap \aleph.
\end{equation*}
To prove the reverse inclusion, it suffices to show that $\sum_{i\in\mN}\DB_{\phi_i}(Y_i, X_i)=\DB_{\bphi}^{n}(\Y, \X)$ for any
\begin{equation*}
\Y\in\{\Z~|~\sum_{i\in\mN}\DB_{\phi_i}(Z_i, X_i)\leq\epsilon\} 
    \cap \aleph\,.
\end{equation*}
From the proof of $i)$, we have $\sum_{i\in\mN}\DB_{\phi_i}(Y_i, X_i)\leq \DB_{\bphi}^{n}(\Y, \X)$ for any random vectors $\X$ and $\Y$. Thus, it remains to show that $\sum_{i\in\mN}\DB_{\phi_i}(Y_i, X_i)\geq \DB_{\bphi}^{n}(\Y, \X)$ for any $\Y\in\aleph$. As the infimum in the univariate BW divergence is attained by the comonotonic coupling, it follows that
\begin{equation*}
\sum_{i\in\mN}\DB_{\phi_i}(Y_i, X_i) = 
    \sum_{i\in\mN}\E[B_{\phi_i}(Y_i, X_i)] = 
\E[B_{\bphi}(\Y, \X)] \geq \DB_{\bphi}^{n}(\Y, \X),
\end{equation*}
where the inequality follows from the definition of the BW divergence. 
\end{proof}

We are now ready to derive lower and upper bounds for DRO problem \eqref{eq:DRO-g-formulation} with uncertainty quantified with a separable Bregman generator when the risk functional $\rho$ is law invariant, subadditive, and comonotonic additive. Recall that a risk functional $\rho$ is subadditive if $\rho(X+Y)\leq\rho(X)+\rho(Y)$ for all rvs $X,Y$ for which the risk functional is well defined. Moreover, a risk functional is comonotonic additive if $\rho(X+Y)=\rho(X)+\rho(Y)$ for all comonotonic rvs $X,Y$ for which the risk functional is well defined.

\begin{proposition} [Worst-case risks for separable Bregman generators]\label{prop:bregman_sum} Suppose $\bphi(\x)=\sum_{i\in\mN}\phi_i(x_i)$ is a separable Bregman generator and $g(\x)=\sum_{i\in\mN}g_i(x_i)$ is a separable aggregation function. Let $\rho$ be a law-invariant risk functional. Then, the following hold for all $\epsilon\geq 0$: 
\begin{enumerate}[label = $\roman*)$]
\item If $\rho$ is subadditive, then
\begin{equation*}
\sup_{\Z\in\bmfB^{n}_{\bphi, \epsilon}(\X)}\rho(g(\Z))
    \leq
\sum_{i\in\mN}\,\sup_{Z_i\in\mfB_{\phi_i, \epsilon}
    (X_i)}\rho \big(g_i(Z_i)\big).
\end{equation*}
\item If $\X$ is comonotonic \footnote{An $n$-dimensional random vector $\Z$ is said to be comonotonic if the rvs $Z_i, Z_j$ are comonotonic for all $i,j\in\mN$ such that $i\neq j$.}, $\rho$ is comonotonic additive, and $g_i$ is non-decreasing for all $i\in\mN$, then
\begin{equation*}
\sum_{i\in\mN}\,\sup_{Z_i\in\mfB_{\phi_i, 
    \frac{\epsilon}{n}}(X_i)}\rho \big(g_i(Z_i)\big)
        \leq
    \sup_{\Z\in\bmfB^{n}_{\bphi, \epsilon}(\X)}\rho(g(\Z)).
\end{equation*}
\item If $\X$ is comonotonic, $g(\x)=\bbeta^T\x$ for some $\bbeta\in\R^n_{\geq}$, and $\rho$ is positive homogeneous, subadditive, and comonotonic additive, then
\begin{equation}\label{eq:bregman_upper_bound_2}
\sum_{i\in\mN}\,\beta_i\;\sup_{Z_i\in\mfB_{\phi_i, 
    \frac{\epsilon}{n}}(X_i)}\rho (Z_i)
\,\leq\,
    \sup_{\Z\in\bmfB^{n}_{\bphi, \epsilon}(\X)}\rho(\bbeta^T\Z)
\,\leq\,
    \sum_{i\in\mN}\,\beta_i\sup_{Z_i\in\mfB_{\phi_i,
\epsilon}(X_i)}\rho (Z_i).
\end{equation}
\end{enumerate}
\end{proposition}

\begin{proof}
We first prove $i)$. Since the Bregman generator $\bphi$ is separable, 
\begin{align*}
\sup_{\Z\in\mfB_{\bphi, \epsilon}^{n}(\X)}\rho\big(g(\Z)\big) 
    &\leq \sup_{\{\Z|\sum_{i\in\mN}
\DB_{\phi_i}(Z_i, X_i)\leq\epsilon\}}
    \rho\big(\sum_{i\in\mN}g_i(Z_i)\big) \\
&\leq \sup_{\{\Z|\sum_{i\in\mN}
    \DB_{\phi_i}(Z_i, X_i)\leq\epsilon\}}
\sum_{i\in\mN}\rho\big(g_i(Z_i)\big) \\
&\leq \sup_{\{\Z|\DB_{\phi_i}(Z_i, X_i)\leq\epsilon 
    \ \text{for all} \ i\in\mN\}}\sum_{i\in\mN} 
\rho\big(g_i(Z_i)\big) \\
&= \sum_{i\in\mN}\sup_{Z_i\in
    \mfB_{\phi_i, \epsilon}(X_i)}\rho\big(g_i(Z_i)\big),
\end{align*}
where the first inequality follows from \Cref{theorem:bregman_sum} $i)$, the second inequality follows from the subadditivity of $\rho$, and the third inequality follows from the fact that $\sum_{i\in\mN}\DB_{\phi_i}(Z_i, X_i)\leq\epsilon$ implies that $\DB_{\phi_i}(Z_i, X_i)\leq\epsilon$ for all $i\in\mN$. 

\smallskip

Second, we prove $ii)$. Suppose that $\X$ is comonotonic and let $\aleph$ be as defined in Theorem \ref{theorem:bregman_sum}.
Then,
\begin{align*}
\sup_{\Z\in\mfB_{\bphi, \epsilon}^{n}(\X)}\rho\big(g(\Z)\big) 
    &\geq \sup_{\Z\in\mfB_{\bphi, \epsilon}^{n}(\X)\cap\aleph}\rho\big(\sum_{i\in\mN}g_i(Z_i)\big) \\
&= \sup_{\{\Z|\sum_{i\in\mN}\DB_{\phi_i}(Z_i, X_i)\leq  
    \epsilon\}\cap\aleph}\rho\big(
\sum_{i\in\mN}g_i(Z_i)\big) \\
&= \sup_{\{\Z|\sum_{i\in\mN}\DB_{\phi_i}(Z_i, X_i)\leq
    \epsilon\}\cap\aleph}
\sum_{i\in\mN}\rho\big(g_i(Z_i)\big) \\
&= \sup_{\{\Z|\sum_{i\in\mN}\DB_{\phi_i}(Z_i, X_i)\leq
    \epsilon\}}\sum_{i\in\mN}\rho\big(g_i(Z_i)\big) \\
&\geq \sup_{\{\Z|\DB_{\phi_i}(Z_i, X_i)\leq\epsilon/n 
    \ \text{for all} \ i\in\mN\}} 
\sum_{i\in\mN}\rho\big(g_i(Z_i)\big) \\
&= \sum_{i\in\mN}\,\sup_{Z_i\in\mfB_{\phi_i, 
    \frac{\epsilon}{n}}(X_i)}
\rho\big(g_i(Z_i)\big), \hphantom{\sum_{i\in\mN}\ \ }
\end{align*}
where the first equality follows from \Cref{theorem:bregman_sum} $ii)$. For the second equality, observe that $\Z$ is comonotonic since $\X$ is comonotonic and $Z_i, X_i$ are comonotonic for all $i\in\mN$. Moreover, as $g_i$ is non-decreasing for all $i\in\mN$, it follows that $(g_1(Z_1), \ldots, g_n(Z_n))$ is comonotonic, so the second equality follows by comonotonic additivity of $\rho$. The third equality follows by law invariance and the second inequality holds since $\DB_{\phi_i}(Z_i, X_i)\leq\frac{\epsilon}{n}$ for all $i\in\mN$ implies that  $\sum_{i\in\mN}\DB_{\phi_i}(Z_i, X_i)\leq\epsilon$.

\smallskip

Finally, the bounds in \eqref{eq:bregman_upper_bound_2} follow from $i)$ and $ii)$ whenever $\rho$ is positive homogeneous. 
\end{proof}

Note that when we use an arbitrary separable Bregman generator, we require different assumptions on the aggregation function $g$ and the reference random vector $\X$. In particular, \Cref{prop:bregman_sum} requires $g$ to be separable and $\X$ must be comonotonic for the lower bound. Comonotonicity of $\X$ implies that there exists a random variable $Y$ and non-decreasing functions $h_i, \, i\in\mN$ such that for all $i\in\mN$, $X_i=h_i(Y)$ (proposition 4.5 \citet{denneberg1994non}). In a setting where $\X$ is a vector of risky asset prices, then $Y$ could represent an underlying market risk factor that affects all of the asset prices in a similar manner. 

From \Cref{prop:bregman_sum}, we see that the lower and upper bounds are with respect to the univariate BW divergence, and the tolerance distances are $\ep/n$ and $\ep$, respectively. Clearly, for dimension $n = 1$ the upper and lower bounds coincide. Moreover, if $n \to +\infty$ and $\phi_i$ is strictly convex for all $i\in\mN$, then the $i$-th uncertainty set in the lower bound converges to the set $\{V~|~V\eqd X_i\}$ for all $i\in\mN$ and hence the lower bound becomes trivial. Thus, the tightness of the bounds is dependent on the ratio between the dimension of the reference random vector and the tolerance distance $\ep$. We refer to \Cref{sec:mult_bound} for a more detailed discussion of the bounds and a comparison to the Wasserstein case. 

\subsection{The Mahalanobis Distance} \label{sec:m_distance}

In \Cref{sec:separable}, we work with an arbitrary separable Bregman generator $\bphi$. Under stronger assumptions on the structure of $\bphi$, we can obtain tighter bounds for the worst-case risk under fewer assumptions on the risk functional $\rho$. One example is the Wasserstein distance discussed in \Cref{sec:Lipschitz}. Here, we present an analogue of \Cref{theorem:Lipschitz} when the Bregman generator is the Mahalanobis distance. The Mahalanobis distance is a popular distance used in classification problems (\citet{mclachlan}) as well as in financial applications, see e.g., \citet{jaimungal2024kullback}.

\begin{definition} [Mahalanobis Distance] \label{definition:m_distance} Let $Q$ be a symmetric, positive-semi-definite  $n\times n$ matrix. Then, the squared Mahalanobis distance between $\x\in\R^n$ and $\y\in\R^n$ is given by $(\x-\y)^TQ(\x-\y)$. 
\end{definition}

The main result of this subsection is for a special case of the Mahalanobis distance when $Q$ is a diagonal matrix. Note that the $\L^2$ norm is a special case of the Mahalanobis distance with $Q$ equal to the identity matrix, i.e., $\bmfB^{n}_{\x^T\x, \epsilon}(\X)=\bU_{\sqrt{\epsilon}}^{n}(\X)$ and $\mfB_{x^2, K^2\epsilon}\big(g(X)\big)=\U_{K\sqrt{\epsilon}}(g(\X))$. Hence, \Cref{theorem:Lipschitz} (with $a=2$ and $\lambda_{\min} = \lambda_{\max} = 1$) is a special case of this result.

\begin{theorem} [Mahalanobis distance]\label{theorem:m_distance} Let $g:\R^n\to\R$ be a $K$-Lipschitz function (wrt the $\L^2$ norm) with decomposition \eqref{eq:decomposition} and $\X\in\L^p_n$ with $g(\X)\in\L^p$. Let $Q$ be a positive-semi-definite diagonal $n\times n$ matrix with at least one non-zero entry and let $\lambda_{\min}\geq 0$ and $\lambda_{\max}>0$ denote the smallest and largest eigenvalues of $Q$, respectively. Then, for any $\epsilon\geq 0$, 
\begin{equation*}
\mfB_{\lambda_{\max}x^2, \norm{\bbeta}_2^2\epsilon}
    \big(g(\X)\big) \subseteq 
g\big(\bmfB^{n}_{\x^TQ\x, \epsilon}(\X)\big)
    \subseteq
\mfB_{\lambda_{\min}x^2, K^2\epsilon}\big(g(\X)\big).
\end{equation*}
\end{theorem}

\begin{proof}
We first prove that $g\big(\bmfB^{n}_{\x^TQ\x, \epsilon}(\X)\big) \subseteq \mfB_{\lambda_{\min}x^2, K^2\epsilon}\big(g(\X)\big)$. Let $\Z\in \bmfB^{n}_{\x^TQ\x, \epsilon}(\X)$. Then, by definition of $\bmfB^{n}_{\x^TQ\x, \epsilon}(\X)$, we have that $\DB^{n}_{\bphi}[F_{\Z}, F_{\X}]\leq\epsilon$. Moreover, the BW divergence from $g(\Z)$ to $g(\X)$ satisfies
\begin{align*}
\DB_{\lambda_{\min}x^2}[F_{g(\Z)}, F_{g(\X)}] &= 
    \inf_{X\eqd g(\X), \ Z\eqd g(\Z)}
\E[\lambda_{\min}(Z-X)^2] \\
&= \inf_{g(\X')\eqd g(\X), \ g(\Z')\eqd g(\Z)}\lambda_{\min}\, 
    \E\big[\abs{g(\Z')-g(\X')}^2\big] \\
&\leq \inf_{g(\X')\eqd g(\X), \ g(\Z')\eqd g(\Z)}  
    \E\big[K^2\lambda_{\min}\norm{\Z'-\X'}_{2}^{2}\big]\\
&\leq \inf_{\X'\eqd \X, \ \Z'\eqd \Z} 
    \E\big[K^2\lambda_{\min}\norm{\Z'-\X'}_{2}^{2}\big] \\
&= \inf_{\X'\eqd \X, \ \Z'\eqd \Z} 
    \E\big[K^2\lambda_{\min}\sum_{i\in\mN}(Z'_i-X'_i)^2\big] \\
&\leq \inf_{\X'\eqd \X, \ \Z'\eqd \Z} 
    \E\big[K^2(\Z'-\X')^TQ(\Z'-\X')\big] \\
&= K^2\,\DB^{n}_{\x^TQ\x}[F_{\Z}, F_{\X}] \\
&\leq K^2\epsilon\,,
\end{align*}
where the first inequality follows by Lipschitz continuity of $g$. Thus, $g(\Z) \in \mfB_{\lambda_{\min}x^2, K^2\epsilon}\big(g(\X)\big)$.

\smallskip

Next, we prove the first set inclusion. Note that since $\lambda_{\max}>0$, $\mfB_{\lambda_{\max}x^2, \norm{\bbeta}_2^2\epsilon}\big(g(\X)\big)=\{V~|~V\eqd g(\X)\}$ whenever $\norm{\bbeta}_2=0$, so we assume wlog that $\norm{\bbeta}_2>0$. Let $Z\in\mfB_{\lambda_{\max}x^2, \norm{\bbeta}_2^2\epsilon}(g(\X))$. Define the $(n-m)$-dimensional vector  
\begin{equation*}
\tilde{\bbeta} := (\sgn(\beta_1)\abs{\beta_1}, \ldots, 
    \sgn(\beta_{n-m})\abs{\beta_{n-m}}),
\end{equation*}
where $\sgn(\cdot)$ is the sign function. Let $\Z=(\Z^{(1)}, \Z^{(2)})\in\R^m\times\R^{n-m}$, where $\Z^{(1)}:=(X_1, \ldots, X_m)$ and 
\begin{equation*}
\Z^{(2)}:=\Big(X_{m+1}+\tfrac{\tilde{\beta}_1
    (Z-g(\X))}{\bbeta^T\tilde{\bbeta}}, \ldots, X_{n} + 
\tfrac{\tilde{\beta}_{n-m}
    (Z-g(\X))}{\bbeta^T\tilde{\bbeta}}\Big).
\end{equation*} 
It suffices to show that $g(\Z)=Z$ and $\DB_{\x^TQ\x}^{n}(F_{\Z}, F_{\X})\leq\epsilon$. The first fact follows from similar calculations as in the proof of \Cref{prop:lin_comp}. For the second claim,
\begin{align*}
\DB_{\x^TQ\x}^{n}(F_{\Z}, F_{\X}) 
    &= \inf_{\X'\eqd\X, \Z'\eqd\Z}
\E\big[(\Z'-\X')^TQ(\Z'-\X')\big] \\
&\leq \E\big[(\Z-\X)^TQ(\Z-\X)\big] \\
&\leq \E\big[\lambda_{\max}
    \sum_{i=1}^{n-m}\frac{\tilde{\beta}_i^2
\big(Z-g(\X)\big)^2}{(\bbeta^T\tilde{\bbeta})^2}\big] \\
&= \frac{\E\big[\lambda_{\max}\big(Z-g(\X)\big)^2\big]}
    {\norm{\bbeta}_2^2} \\
&\leq \frac{\norm{\bbeta}_2^2\epsilon}{\norm{\bbeta}_2^2} \\
&=\epsilon,
\end{align*}
where the second equality follows from similar calculations as in the proof of \Cref{prop:lin_comp}. 
\end{proof}

The following corollary follows from \Cref{theorem:m_distance}.

\begin{corollary} \label{cor:m_distance} Suppose the assumptions of \Cref{theorem:m_distance} are satisfied. Then, for any $\epsilon \geq 0$ and law-invariant risk functional $\rho$,
\begin{equation*} 
\sup_{Y\in \mfB_{\lambda_{\max}x^2, 
    \norm{\bbeta}_2^2\epsilon}(g(\X))}\rho(Y) \leq 
\sup_{\Y\in \bmfB^{n}_{\x^TQ\x, \epsilon}(\X)} 
    \rho\big(g(\Y)\big) \leq 
\sup_{Y\in \mfB_{\lambda_{\min}x^2, 
    K^2\epsilon}(g(\X))}\rho(Y).
\end{equation*}
\end{corollary}

Since the (univariate) Bregman divergence satisfies $B_{c\phi}\big(\z_1, \z_2\big)  = c\, B_{\phi}\big(\z_1, \z_2\big) $ for $c \ge 0$, the lower bound can be rewritten as a DRO problem over a Wasserstein uncertainty set as follows
\begin{equation*}
    \sup_{Y\in \mfB_{\lambda_{\max}x^2, 
    \norm{\bbeta}_2^2\epsilon}(g(\X))}\rho(Y)
    =
    \sup_{Y\in \mfB_{x^2, 
    \frac{\norm{\bbeta}_2^2\epsilon}{\lambda_{\max}}}(g(\X))}\rho(Y)
    =
    \sup_{Y\in \U_{\ep'}(g(\X))}\rho(Y)\,,
\end{equation*}
where $\ep' := \norm{\bbeta}_2\sqrt{\frac{\epsilon}{\lambda_{\max}}}$. 
Similarly the upper bound can be recast as a Wasserstein DRO problem over the uncertainty set $\U_{\ep^\dagger}(g(\X))$ with $\ep^\dagger := K\sqrt{\frac{\epsilon}{\lambda_{\min}}}$.

We further see that the tolerance distance for the lower bound is proportional to the norm of $\bbeta$ (i.e., the linear term of the aggregation function), while for the upper bound it is proportional to the Lipschitz constant of $g$. Thus, the bounds are tighter, if the linear term in $g$ dominates (i.e., $K^2 - \norm{\bbeta}_2^2$ is small) and the eigenvalues of $Q$ are close to each other. Furthermore, when $g$ is linear and $\lambda_{\min}=\lambda_{\max}$, the inequalities become equalities. This is not surprising since the case when $\lambda_{\min}=\lambda_{\max}$ corresponds to a scaled version of the squared Wasserstein distance. 

\subsection{Composable Bregman Generators} \label{sec:composable}

Here, we consider composable Bregman generators $\bphi:\mathcal{X}\to\R$ of the form $\bphi(\x)=\phi(g(\x))$, where $g:\mathcal{X}\to\R$ and $\phi:\R\to\R$ are Bregman generators and $\phi$ is non-decreasing. Non-decreasingness of $\phi$ (combined with convexity of $g$), ensures that $\bphi$ is convex.

Recall that for a DRO problem $\sup_{\X\in \U}\rho\big(g(\X)\big)$ with some uncertainty set $\U$ to be well defined, the aggregation function $g\colon \R^n \to \R$, mapping input to output, is often convex and $\U$ is a convex set. Thus, composable Bregman generators allow us to tailor the uncertainty set to the aggregation function at hand. Indeed, for a univariate and non-decreasing Bregman generator $\phi$, the uncertainty set induced by the composable BW divergence, i.e., $\bmfB^{n}_{\phi\circ g, \epsilon}(\X)$, is quantified via the Bregman divergence as a composition of the aggregation function and $\phi$. Thus, the aggregation function plays a pivotal role in defining the uncertainty set. 

For composable Bregman generators, we again bound the risk-aware DRO problem with multivariate BW uncertainty with a risk-aware DRO problem with univariate BW uncertainty. Depending on the choice of law-invariant risk functional, the univariate DRO problem can be solved analytically. We refer to \Cref{sec:univ_bound} where we solve the univariate BW DRO problem when $\rho$ is a signed Choquet integral.
\begin{theorem} [Composable Bregman generators]\label{theorem:bregman_upper_bound}
Let $\bphi:\mathcal{X}\to\R$ be a Bregman generator of the form $\bphi(\x)=\phi(g(\x))$, where $g:\mathcal{X}\to\R$, $\phi:\R\to\R$ are Bregman generators and $\phi$ is non-decreasing. Then, for any law-invariant risk functional and  $\epsilon\geq 0$, 
\begin{equation*}
\sup_{\Z\in\bmfB^{n}_{\bphi, \epsilon}(\X)}\rho(g(\Z))
    \,\leq\,
\sup_{Z\in\mfB_{\phi, \epsilon}(g(\X))}\rho(Z).
\end{equation*}
\end{theorem}

\begin{proof}
It suffices to show that $\Z\in\bmfB^{n}_{\bphi, \epsilon}(\X)$ implies that $ g(\Z)\in\mfB_{\phi, \epsilon}(g(\X))$. For any random vector $\Z$,
\begin{align*}
\DB^{n}_{\bphi}[F_{\Z},F_{\X}] &= 
    \inf_{\Z'\eqd\Z, \X'\eqd\X}
\E[\bphi(\Z')-\bphi(\X')-\nabla\bphi(\X')\cdot(\Z'-\X')] \\
&= \inf_{\Z'\eqd\Z, \X'\eqd\X}\E[\phi(g(\Z'))-\phi(g(\X')) 
    -\phi'(g(\X'))\nabla g(\X')\cdot(\Z'-\X')] \\
&= \inf_{\Z'\eqd\Z, \X'\eqd\X}\E\big[\phi(g(\Z'))-\phi(g(\X')) 
    -\phi'(g(\X'))\nabla g(\X')\cdot(\Z'-\X') \\
&\qquad\qquad\qquad\quad -\phi'(g(\X'))(g(\Z')-g(\X')) + 
    \phi'(g(\X'))(g(\Z')-g(\X'))\big] \\
&= \inf_{\Z'\eqd\Z, \X'\eqd\X}\E[B_{\phi}(g(\Z'), g(\X')) + 
    \phi'(g(\X'))B_g(\Z',\X')] \\
&\geq \DB_{\phi}[F_{g(\Z)}, F_{g(\X)}] + 
    \inf_{\Z'\eqd\Z, \X'\eqd\X}\E[\phi'(g(\X'))B_g(\Z',\X')] \\
&\geq \DB_{\phi}[F_{g(\Z)}, F_{g(\X)}], 
\end{align*}
where the last inequality follows from the fact that $\phi'(g(\X'))$ is non-negative (since $\phi$ is non-decreasing) and the Bregman divergence is non-negative. Thus, for any $\Z\in\bmfB^{n}_{\bphi, \epsilon}(\X)$, we have that $\DB_{\phi}[F_{g(\Z)}, F_{g(\X)}]\leq\DB^{n}_{\bphi}[F_{\Z},F_{\X}]\leq\epsilon$. 
\end{proof}

From the proof, we observe that the quality of the bound depends on the random variable $\phi'(g(\X'))B_g(\Z',\X')$. In particular, the bound is better when the expected value of $\phi'(g(\X'))B_g(\Z',\X')$ is small. Moreover, if $B_g(\Z',\X') = 0$, then equality in \Cref{theorem:bregman_upper_bound} holds; however this is a strong requirement. 

\section{Explicit Bounds for Signed Choquet Integrals} \label{sec:application}

Next, we derive explicit bounds for the worst-case risk when the risk functional belongs to the class of signed Choquet integrals. In \Cref{sec:univ_bound}, we characterize the distribution that attains the largest risk over a univariate BW uncertainty set. Then, in \Cref{sec:mult_bound}, we apply these results to derive formulae for the bounds established in Sections \ref{sec:Wasserstein} and \ref{sec:bw}. We first recall the definition of a signed Choquet integral, which was introduced by \citet{choquet}. 

\begin{definition} [Signed Choquet Integrals] \label{definition:choquet} A signed Choquet integral, denoted $I_h$, is a mapping from $\L^{p}$ to $\mathbb{R}$ given by 
\begin{equation} \label{eq:choquet}
I_h(X)=\int_{-\infty}^{0}\big[h(\mathbb{P}(X\geq x))-h(1)\big]\, dx+\int_{0}^{\infty}h(\mathbb{P}(X\geq x))\, dx,
\end{equation}
where $h: [0,1]\to\R$ is non-constant, has bounded variation, and satisfies $h(0)=0$. The function $h$ is called the distortion function of $I_h$. 
\end{definition}

Note that the integrals in \Cref{eq:choquet} may be infinite, in which case $I_h(X)$ may not be well defined. In this paper, we assume that $X$ is such that $I_h(X)$ is finite and thus well defined, and refer to \citet{wang} for an extensive discussion on signed Choquet integrals. When the distortion function $h$ is non-decreasing and satisfies $h(1)=1$, then $I_h(X)$ is a distortion risk measure. The class of distortion risk measures includes many commonly used risk measures, including the Value-at-Risk (VaR) and the ES. However, several important risk measures such as the Interquartile Range, Mean-Median Difference, and Gini Deviation belong to the class of signed Choquet integrals but are not distortion risk measures. Furthermore, inverse S-shaped distortions, which are popular in economics (\citet{tversky}), belong to the class of signed Choquet integrals. 

If the distortion function $h$
is absolutely continuous, then by lemma 3 in \citet{wang}, the signed Choquet integral has representation
\begin{equation} \label{eq:quantile_rep}
I_h(X)=\int_{0}^{1}\gamma (u)F_X^{-1}(u)\,du\,,
\end{equation}
where $\gamma : [0,1]\to \R$ is called a distortion weight function and defined by $\gamma(u) := \frac{d^-}{dx} h(x)|_{x = 1-u}$, where $\frac{d^-}{dx}$ denotes the left derivative. Throughout the exposition, we only work with signed Choquet integrals with representation \eqref{eq:quantile_rep}. Moreover, as signed Choquet integrals are law invariant, we use the notation $\tilde{I}_h(F_X^{-1}):=I_h(X)$ whenever $X$ has quantile function $F_X^{-1}$. Additionally, for the class of distortion risk measures, the  distortion weight function $\gamma$ is non-negative (since $h$ is non-decreasing) and $\int_{0}^{1}\gamma (u)\,du=h(1)-h(0)=1$, thus $\gamma$ is a density on $[0,1]$.

\subsection{Worst-case Quantile Function for Univariate Risks} \label{sec:univ_bound}

The largest signed Choquet integral under univariate BW uncertainty is given by

\begin{equation} \label{opt:WC-RMtoto}
\max_{\{G^{-1}~|~G\in\mathcal{M}_p(\R)\}}\tilde{I}_h(G^{-1}), \quad \text{subject to} \quad \DB_\phi (G^{-1}, F^{-1})\leq\epsilon\,,
\end{equation}
where $\DB_\phi (G^{-1}, F^{-1}):=\DB_\phi [G, F]$.
We call the quantile function attaining the maximum the worst-case quantile function. The power $p$ of the feasible set $\mathcal{M}_p(\R)$ must be chosen such that the divergence constraint is well defined. For example, for the squared $2$-Wasserstein distance, $p=2$ is suitable. 

For $\phi(x)=x^2$, optimization problem \eqref{opt:WC-RMtoto} coincides with the upper and lower bounds given by \Cref{prop:Lipschitz} (when $a=p=2$) for an appropriately chosen tolerance distance. Moreover, for a separable Bregman generator $\bphi(\x)=\sum_{k\in\mN}\phi_k(x_k)$, we can compute the upper and lower bounds in Proposition \ref{prop:bregman_sum} by solving the optimization problem \eqref{opt:WC-RMtoto} for each $\phi_k$ separately. 

While optimization problem \eqref{opt:WC-RMtoto} with the Wasserstein distance and distortion risk measures has been studied in \citet{liu2022inf}, little has been done with BW uncertainty. An exception is \citet{pesenti}, which studies optimization problem \eqref{opt:WC-RMtoto} for non-decreasing and strictly concave distortion functions $h$ (i.e., $\gamma$ non-negative and strictly decreasing). We generalize their result in two directions. First, we consider concave distortion functions $h$ that are not necessarily non-decreasing --- corresponding to the class of subadditive signed Choquet integrals. Second, we consider non-decreasing $h$ that are not necessarily concave --- corresponding to monotone signed Choquet integrals.

To solve optimization problem \eqref{opt:WC-RMtoto} for non-concave $h$, we recall the isotonic projection, which is the continuous analogue to the isotonic regression; see, e.g., \citet{brunk, barlow, barlow2}. The use of isotonic projections to solve problems related to distortion risk measures can for example be found in \citet{pesenti2} and \citet{bernard}. For this we write $\L^2((0,1)):=\{\ell \colon (0,1) \to \R ~|~ \int_0^1 \ell^2(u) \, du < \infty\}$ to denote the set of all square-integrable functions defined on $(0,1)$. Moreover, for any $\ell\in\L^2((0,1))$, we denote by $\norm{\ell}_2:=\sqrt{\int_{0}^{1}\ell^2(u)\,du}$ its $\mathcal{L}^2$ norm.

\begin{definition} [Isotonic Projection] \label{definition:isotonic_proj} The isotonic projection of a function $\ell\in\L^2((0,1))$, denoted $\iso{\ell}$, is
\begin{equation*}
\iso{\ell}:=\argmin_{j\in\K}\int_{0}^{1}\big(j(u)-\ell(u)\big)^2du,
\end{equation*}
where $\K\subseteq \mathcal{L}^2((0,1))$ is the set of all left-continuous and non-decreasing square-integrable functions. 
\end{definition}
Note that $\K$ can be identified as the set of quantile functions. Furthermore, the isotonic projection preserves the ordering of functions in $\L^2((0,1))$, discussed next. 

\begin{proposition} [Ordering property]\label{prop:iso_ordering}
Let $\ell_1, \ell_2\in\L^2((0,1))$. If $\ell_2(s)\leq \ell_1(s)$ for all $s\in (0,1)$, then $\iso{\ell_2}(s)\leq\iso{\ell_1}(s)$. This result also holds when all inequalities are replaced with strict inequalities. 
\end{proposition}

\begin{proof}
See Appendix \ref{appendix:b}.
\end{proof}

\begin{theorem} [Worst-case quantile function]\label{theorem:wc_signed_choquet_int} 
Let $I_h$ be a signed Choquet integral with representation \eqref{eq:quantile_rep} and distortion weight function $\gamma$. Suppose $\phi$ is a strictly convex Bregman generator satisfying $\lim_{x\to\infty}\phi'(x)=\infty$ and $\lim_{x\to -\infty}\phi'(x)=-\infty$. Assume at least one of the following holds:

\begin{enumerate}[label = $\roman*)$]
    \item $\gamma$ is non-decreasing
    \item $\gamma$ is non-negative and $\int_0^1\phi'\big(F^{-1}(u)\big)^2 + \gamma(u)^2 \,du< \infty$.
\end{enumerate}
Then, the worst-case quantile function of \eqref{opt:WC-RMtoto}, if it exists, is unique and given by
        \begin{equation}\label{eq:solution1}
        G^{-1}_{\lambda^*}(u)
        := 
        \left(\phi'\right)^{-1}
        \left(\left(
        \phi'\big(F^{-1}(u)\big) + \frac{1}{\lambda^*}\gamma(u)
        \right)^{\uparrow}\right)
        \,,
        \end{equation}
    where $\lambda^{*}>0$ is the smallest solution to $\DB_\phi\big(G_{\lambda}^{-1}, \, F^{-1}\big) = \epsilon$ in case $i)$ and the unique solution to $\DB_\phi\big(G_{\lambda}^{-1}, \, F^{-1}\big) = \epsilon$ in case $ii)$.
\end{theorem}
We note that in case $i)$, the function $\phi'\big(F^{-1}(u)\big) + \frac{1}{\lambda^*}\gamma(u)$ is non-decreasing and equal to its isotonic projection.

\begin{proof}
Suppose either $i)$ or $ii)$ be satisfied. We first show that the solution (if it exists) is of the form \eqref{eq:solution1}. Let the worst-case quantile function have a BW constraint of $\ep_0 \in [0, \ep]$. Then the Lagrangian for the constrained optimization problem \eqref{opt:WC-RMtoto} is 
\begin{equation*}
\mathcal{L}(G^{-1}, \lambda)=\lambda\int_{0}^{1}\phi\left(G^{-1}(u)\right)-\Big(\frac{1}{\lambda}\gamma(u)+\phi'\left(F^{-1}(u)\right)\Big)G^{-1}(u)du+k(\lambda),
\end{equation*}
where $k(\lambda)$ is a function that does not depend on $G^{-1}$.
Since $\phi$ is strictly convex, by \citet{barlow} theorem 3.1, the argmin of the Lagrangian over the set of quantile functions $\K$ is uniquely attained at $\tilde{G}^{-1}_{\lambda}(u):=\left(\phi'\right)^{-1}\left(\left(\phi'\big(F^{-1}(u)\big)+\frac{1}{\lambda}\gamma(u)\right)^{\uparrow}\right)$.

\smallskip

For case $i)$, since $\gamma$ is non-decreasing, $F^{-1}(u)$ is non-decreasing, and $\phi'$ is strictly increasing, the function $\tilde{G}^{-1}_{\lambda}(u)=(\phi')^{-1}(\phi'(F^{-1}(u))+\frac{1}{\lambda}\gamma(u))$ is a non-decreasing function of $u$ (and hence equal to its isotonic projection). Next, we prove existence of the Lagrange parameter for case $i)$. For this we first show that the signed Choquet integral is strictly decreasing with respect to the Lagrange parameter, i.e., $0<\lambda_1<\lambda_2 $ implies that $\tilde{I}_h(\tilde{G}_{\lambda_1}^{-1})>\tilde{I}_h(\tilde{G}_{\lambda_2}^{-1})$. Note that $0<\lambda_1<\lambda_2$ implies that
\begin{equation}\label{eq:G}
\begin{cases} \tilde{G}_{\lambda_1}^{-1}(u)> \tilde{G}_{\lambda_2}^{-1}(u) & \text{if} \quad \gamma(u)> 0 \\[0.5em] \tilde{G}_{\lambda_1}^{-1}(u)< \tilde{G}_{\lambda_2}^{-1}(u) & \text{if} \quad \gamma(u)<0 \,.
\end{cases} 
\end{equation}
Since $\gamma$ is non-decreasing, there exists $y\in [0,1]$ such that $\gamma(u)<0$ on $[0,y)$ and $\gamma(u)\geq 0$ on $(y,1]$. Therefore, we can rewrite the signed Choquet integral as 
\begin{equation*}
\tilde{I}_h(\tilde{G}_{\lambda}^{-1})=\int_{0}^{y}\gamma(u)\tilde{G}_{\lambda}^{-1}(u)du+\int_{y}^{1}\gamma(u)\tilde{G}_{\lambda}^{-1}(u)du,
\end{equation*}
and it follows from \eqref{eq:G} that $\tilde{I}_h(\tilde{G}_{\lambda}^{-1})$ is a strictly decreasing function of $\lambda$. Thus, the optimal Lagrange parameter is the smallest $\lambda>0$ such that the constraint is satisfied. Finally, as $\phi$ is differentiable and strictly convex, it is continuously differentiable and therefore $\DB_\phi\big(\tilde{G}_{\lambda}^{-1}, \, F^{-1}\big)$ is continuous in $\lambda$. Furthermore, $\lim_{\lambda\to 0}\DB_\phi\big(\tilde{G}_{\lambda}^{-1}, \, F^{-1}\big)=\infty$
and $\lim_{\lambda\to \infty}\DB_\phi\big(\tilde{G}_{\lambda}^{-1}, \, F^{-1}\big)=0$. Thus, a unique smallest solution to $\DB_\phi\big(\tilde{G}_{\lambda}^{-1}, \, F^{-1}\big) = \epsilon$ exists for any $\epsilon>0$.

\smallskip

Next, we prove existence of the Lagrange parameter for case $ii)$. Again, we first show that the signed Choquet integral is strictly decreasing in $\lambda$. Let $0<\lambda_1<\lambda_2$. Since $\gamma$ is non-negative, by \Cref{prop:iso_ordering}, 
\begin{equation*}
    \tilde{G}^{-1}_{\lambda_1}(u)=\left(\phi'\right)^{-1}\Big(\Big(\phi'\big(F^{-1}(u)\big)+\frac{1}{\lambda_1}\gamma(u)\Big)^{\uparrow}\;\Big)
    \geq \left(\phi'\right)^{-1}\Big(\Big(\phi'\big(F^{-1}(u)\big)+\frac{1}{\lambda_2}\gamma(u)\Big)^{\uparrow}\;\Big)=\tilde{G}^{-1}_{\lambda_2}(u),
\end{equation*}
with strict inequality for any $u\in(0,1)$ such that $\gamma(u)\neq 0$. Thus, $\tilde{I}_h(\tilde{G}_{\lambda}^{-1})$ is a strictly decreasing function of $\lambda$ and the optimal Lagrange parameter is the smallest $\lambda>0$ such that the constraint is satisfied. By the same arguments as in case $i)$, a $\lambda^*$ satisfying the constraint always exists. To show uniqueness, it suffices to show that $\DB_\phi\big(\tilde{G}_{\lambda}^{-1}, \, F^{-1}\big)$ is a strictly decreasing function of $\lambda$. For $0<\lambda_1<\lambda_2$, we have
\begin{equation*}
\DB_\phi\big(\tilde{G}_{\lambda_1}^{-1}, \, F^{-1}\big)-\DB_\phi\big(\tilde{G}_{\lambda_2}^{-1}, \, F^{-1}\big)=\int_{S}\big(\tilde{G}_{\lambda_1}^{-1}(u)-\tilde{G}_{\lambda_2}^{-1}(u)\big)\Big(\tfrac{\phi(\tilde{G}_{\lambda_1}^{-1}(u))-\phi(\tilde{G}_{\lambda_2}^{-1}(u))}{\tilde{G}_{\lambda_1}^{-1}(u)-\tilde{G}_{\lambda_2}^{-1}(u)}-\phi'\big(F^{-1}(u)\big)\Big)du,
\end{equation*}
where $S=\{u\in(0,1)~|~\gamma(u)\neq 0\}$. 
Since $\gamma$ is strictly positive on $S$, the first term in the integrand is strictly positive on $S$. For the second term in the integrand,
\begin{align*}
\frac{\phi\big(\tilde{G}_{\lambda_1}^{-1}(u)\big)-\phi\big(\tilde{G}_{\lambda_2}^{-1}(u)\big)}{\tilde{G}_{\lambda_1}^{-1}(u)-\tilde{G}_{\lambda_2}^{-1}(u)}-\phi'\big(F^{-1}(u)\big) &> \phi'\big(\tilde{G}_{\lambda_2}^{-1}(u)\big)-\phi'\big(F^{-1}(u)\big) \\
&= \Big(\phi'\big(F^{-1}(u)\big) + \frac{1}{\lambda_2}\gamma(u)
        \Big)^{\uparrow}-\phi'\big(F^{-1}(u)\big) \\
&\geq \big(\phi'\big(F^{-1}(u)\big)\big)^{\uparrow}-\phi'\big(F^{-1}(u)\big) \\
&= 0,
\end{align*}
where the first inequality follows from the strict convexity of $\phi$, the equality follows from \eqref{eq:solution1}, the second inequality follows from \Cref{prop:iso_ordering}, and the last equality follows from the fact that $\phi'\big(F^{-1}(u)\big)$ is non-decreasing. 
\end{proof}

\subsection{Bounds for Multivariate Risks} \label{sec:mult_bound}

In this subsection, we calculate the bounds for worst-case signed Choquet integrals under the settings considered in Sections \ref{sec:Wasserstein} and \ref{sec:bw}. In particular, we derive explicit bounds for the multivariate Wasserstein setup (\Cref{prop:Lipschitz}), BW divergences with separable Bregman generators (\Cref{prop:bregman_sum}), and BW divergences with the Mahalnobis distance generator (\Cref{cor:m_distance}).

We start with uncertainty sets defined via the multivariate Wasserstein distance.

\begin{proposition}[Multivariate Wasserstein distance]\label{prop:wc_signed_choquet_int} Suppose that $I_h$ satisfies the conditions of \Cref{theorem:wc_signed_choquet_int}. Further, assume that the aggregation function $g:\R^n\to\R$ is $K$-Lipschitz with decomposition \eqref{eq:decomposition} and the distortion weight function $\gamma$ satisfies $0<\norm{\gamma}_2<\infty$. Let $\X\in\L^2_n$ such that $g(\X)\in\L^2$ and $\epsilon>0$.  
\begin{enumerate}[label = $\roman*)$]
\item If $\gamma$ non-decreasing, then 
\begin{equation*}
I_h(g(\X))+\norm{\bbeta}_2\norm{\gamma}_2\epsilon\leq\sup_{\Y\in \bU^{n}_{\epsilon}(\X)}I_h\big(g(\Y)\big)\leq I_h(g(\X))+K\norm{\gamma}_2\epsilon. 
\end{equation*}
\item If $\gamma$ is non-negative, then 
\begin{equation} \label{eq:wc_bound_2} \int_{0}^{1}\gamma(u)\iso{\Big(F_{g(\X)}^{-1}(u)+\frac{1}{2\underline{\lambda}}\gamma(u)\Big)}\leq\sup_{\Y\in \bU^{n}_{\epsilon}(\X)}I_h\big(g(\Y)\big)\leq\int_{0}^{1}\gamma(u)\iso{\Big(F_{g(\X)}^{-1}(u)+\frac{1}{2\overbar{\lambda}}\gamma(u)\Big)}du, 
\end{equation}
where $\underline{\lambda}$ is the unique positive solution to 
\begin{equation*}
 \int_0^1 \left(F_{g(\X)}^{-1}(u) - \iso{\Big(F_{g(\X)}^{-1}(u)+\frac{1}{2\underline{\lambda}}\gamma(u)\Big)}\right)^2 du = \norm{\bbeta}_2^2 \, \ep^2  \,,
\end{equation*}
and $\overbar{\lambda}$ is the unique positive solution to
\begin{equation} \label{eq:constraint2}
 \int_0^1 \left(F_{g(\X)}^{-1}(u) - \iso{\Big(F_{g(\X)}^{-1}(u)+\frac{1}{2\overbar{\lambda}}\gamma(u)\Big)}\right)^2 du = K^2 \, \ep^2  \,.
\end{equation}
\end{enumerate}
\end{proposition}

\begin{proof}
For $i)$ we apply \Cref{prop:Lipschitz}, $i)$ to obtain
\begin{equation}\label{eq:proof-bound-multi-wasser}
    \sup_{Y\in \U_{\norm{\bbeta}_2\epsilon}(g(\X))}I_h(Y)\leq\sup_{\Y\in \bU^{n}_{\epsilon}(\X)}I_h\big(g(\Y)\big)\leq
  \sup_{Y\in \U_{K\epsilon}(g(\X))}I_h(Y)\,.  
\end{equation}
It suffices to explicitly compute the bounds in \eqref{eq:proof-bound-multi-wasser}. We start with the upper bound. 

\smallskip

By \Cref{theorem:wc_signed_choquet_int} $i)$ with $\phi(x)=x^2$, the quantile function attaining the upper bound in \eqref{eq:proof-bound-multi-wasser} is of the form \eqref{eq:solution1}. Plugging $\phi(x)=x^2$ and $F^{-1}(u)=F^{-1}_{g(\X)}(u)$ into \eqref{eq:solution1} yields 
\begin{equation} \label{eq:quantile1}
G_{\lambda}^{-1}(u)=F_{g(\X)}^{-1}(u)+\frac{1}{2\lambda}\gamma(u).
\end{equation}
Moreover, the optimal $\lambda\ge 0$ is the smallest positive solution to 
\begin{equation*}
   W(G_{\lambda}, F_{g(\X)})^2
    =
    \frac{1}{4\lambda^2}\int_0^1 \gamma(u)^2 \, du = K^2  \, \ep^2.
\end{equation*}
Solving for the positive root of $\lambda$, we obtain that $\lambda^{*}=\frac{\norm{\gamma}_2}{2K\epsilon}>0$. Plugging $\lambda^{*}$ into \Cref{eq:quantile1} yields
\begin{equation} \label{eq:wc_univariate}
    \sup_{Y\in \U_{K\epsilon}(g(\X))}I_h(Y)
    =
    I_h\big(g(\X)\big) + K\norm{\gamma}_2\epsilon\,.
\end{equation}
The lower bound follows by the same arguments, with $\norm{\bbeta}_2$ instead of $K$.

\smallskip

For $ii)$, similarly to part $i)$, we apply \Cref{prop:Lipschitz}, $i)$ to obtain \eqref{eq:proof-bound-multi-wasser}. We again start with computing the upper bound. As $g(\X)\in\L^2$ and $\gamma$ is square integrable, the integrability assumption of  \Cref{theorem:wc_signed_choquet_int} part $ii)$ holds. Thus, by \Cref{theorem:wc_signed_choquet_int} part $ii)$, the worst-case quantile function of the upper bound in \eqref{eq:proof-bound-multi-wasser} is 
\begin{align*} 
G_{\lambda}^{-1}(u) &=\frac{1}{2}\iso{\big(2F_{g(\X)}^{-1}(u)+\frac{1}{\lambda}\gamma(u)\big)} =\iso{\big(F_{g(\X)}^{-1}(u)+\frac{1}{2\lambda}\gamma(u)\big)},
\end{align*}
where the second equality follows from  \Cref{prop:iso_properties}, part $i)$, see Appendix \ref{appendix:b}. Again by \Cref{theorem:wc_signed_choquet_int}, we obtain the equation for the (unique) Lagrange multiplier $\overbar{\lambda}$. Calculating the signed Choquet integral with the worst-case quantile function $G_{\overbar{\lambda}}^{-1}$ yields \Cref{eq:wc_bound_2} and the lower bound follows by similar arguments. 
\end{proof}

Part $i)$ above generalizes some well-known results in the literature. For example, when $g$ is linear, then we recover proposition 2 $iv)$ from \citet{wozabal2014robustifying}. Furthermore, when $I_h$ is the expected value (i.e., $\norm{\gamma}_2=1$), we recover theorem 5 from \citet{kuhn} and theorem 6.3 from \citet{esfahani2018data}. Additionally, \eqref{eq:wc_univariate} is a special case of proposition 
4 $i)$ from \citet{liu2022inf}. Using their result, we can generalize $i)$ to the $p$-Wasserstein distance for any $p\geq 1$ provided that $\X\in\L^p_n$ and $g(\X)\in\L^p$. The assumption that the $\L^2$ norm is used in the multivariate Wasserstein distance can also be relaxed. In other words, we can rewrite $i)$ more generally as
\begin{equation} \label{eq:wass_bound_p}
I_h(g(\X))+\norm{\bbeta}_b\norm{\gamma}_q\epsilon\leq\sup_{\Y\in \bU^{n}_{\epsilon}(\X)}I_h\big(g(\Y)\big)\leq I_h(g(\X))+K\norm{\gamma}_q\epsilon,
\end{equation}
where $\frac{1}{p}+\frac{1}{q}=1$ and $\frac{1}{a}+\frac{1}{b}=1$. Here, the $\L^a$ norm is used in the multivariate Wasserstein distance and $K$ is the Lipschitz constant wrt the same norm.

We see that if $\gamma$ is non-decreasing, then the worst-case risk is bounded by the signed Choquet integral of the reference rv plus a positive penalty term which depends on the choice of risk measure. The penalty term in both the upper and lower bounds depend on the level of uncertainty $\epsilon$ in an intuitive way. In particular, increasing the level of uncertainty increases the worst-case risk linearly in $\epsilon$. Additionally, the penalty term of the upper bound depends on the Lipschitz constant of the aggregation function while the lower bound is controlled by the linear components of $g$ (if they exist). In the case where $g$ has no linear components, the worst-case risk is (trivially) lower bounded by the risk of the reference. In contrast, if $g$ is linear, then the lower bound is equal to the upper bound. In other words, the quality of the bounds are controlled by $(K-\norm{\bbeta}_2)\norm{\gamma}_2\epsilon$, where  $K-\norm{\bbeta}_2$ quantifies the ``closeness" of $g$ to a linear function. Thus, the smaller the non-linear components of $g$ (relative to the linear ones), the tighter the bounds. Note that the bounds in \Cref{prop:wc_signed_choquet_int} also hold for $g$ locally Lipschitz as long as $\X$ satisfies the conditions of \Cref{prop:loc_lip}. 

Next, we consider BW divergences with separable Bregman generators and linear aggregation functions $g(\x)=\bbeta^T\x$, $\bbeta\in\R^n_{\geq}$. Linear aggregation function are traditionally studied in DRO problems with Wasserstein uncertainty, indicatively see \citet{wozabal2014robustifying} and \citet{mao}.

\begin{proposition} [Separable Bregman Generators]
\label{prop:sep_Bregman}
Let $I_h$ be a signed Choquet integral with non-decreasing distortion weight function $\gamma$. For $k\in\mN$, let $\phi_k:\mathcal{X}\to\R$ be strictly convex Bregman generators with $\lim_{x\to\infty}\phi_k'(x)=\infty$ and $\lim_{x\to -\infty}\phi_k'(x)=-\infty$ and define the separable Bregman generator $\bphi(\x)=\sum_{k\in\mN}\phi_k(x_k)$. Let $g(\x)=\bbeta^T\x$ for some $\bbeta\in\R^n_{\geq}$. Then, for any $\epsilon>0$ and $\X$ comonotonic,
\begin{align*}
\sum_{i\in\mN}\beta_i\int_{0}^{1}\gamma(u)(\phi_i')^{-1}\Big(\phi_i'\big(F_{X_i}^{-1}(u)\big)+\tfrac{1}{\overbar{\lambda}_i}\gamma(u)\Big)\,du &\leq \sup_{\Z\in\bmfB^{n}_{\bphi, \epsilon}(\X)}I_h(\bbeta^T\Z) \\
&\leq\sum_{i\in\mN}\beta_i\int_{0}^{1}\gamma(u)(\phi_i')^{-1}\Big(\phi_i'\big(F_{X_i}^{-1}(u)\big)+\tfrac{1}{\underline{\lambda}_i}\gamma(u)\Big)\,du,
\end{align*}
where $G_{\lambda_i}^{-1}(u):=(\phi_i')^{-1}\Big(\phi_i'\big(F_{X_i}^{-1}(u)\big)+\frac{1}{\lambda_i}\gamma(u)\Big)$, $\overbar{\lambda}_i$ is the smallest solution to $\DB(G_{\overbar{\lambda}_i}^{-1}, F_{X_i}^{-1})=\epsilon$, and $\underline{\lambda}_i$ is the smallest solution to $\DB(G_{\underline{\lambda}_i}^{-1}, F_{X_i}^{-1})=\frac{\epsilon}{n}$. 

\smallskip

Furthermore, if each $(\phi_i')^{-1}$ is $L_i$ Lipschitz, then the difference between the upper and lower bounds is at most 
\begin{equation*}
\sum_{i\in\mN}\beta_i L_i(\tfrac{1}{\overbar{\lambda}_i}-\tfrac{1}{\underline{\lambda}_i})\norm{\gamma}_2^2.
\end{equation*}
\end{proposition}

\begin{proof}
The proof of the bounds is omitted as it follows similar steps as the proof of \Cref{prop:wc_signed_choquet_int}, with the difference that we apply \Cref{prop:bregman_sum} instead of \Cref{prop:Lipschitz} in the first step.

\smallskip

For the distance between the upper and lower bounds note that
\begin{align*}
&\sum_{i\in\mN}\beta_i\int_{0}^{1}\gamma(u)\Big[(\phi_i')^{-1}\Big(\phi_i'\big(F_{X_i}^{-1}(u)\big)+\tfrac{1}{\overbar{\lambda}_i}\gamma(u)\Big)-(\phi_i')^{-1}\Big(\phi_i'\big(F_{X_i}^{-1}(u)\big)+\tfrac{1}{\underline{\lambda}_i}\gamma(u)\Big)\Big]\,du \\
\leq &\sum_{i\in\mN}\beta_i\int_{0}^{1}\gamma^2(u)L_i(\tfrac{1}{\overbar{\lambda}_i}-\tfrac{1}{\underline{\lambda}_i})\,du \\
= &\sum_{i\in\mN}\beta_i L_i(\tfrac{1}{\overbar{\lambda}_i}-\tfrac{1}{\underline{\lambda}_i})\norm{\gamma}_2^2,
\end{align*}
where the inequality follows by Lipschitz continuity of $(\phi'_i)^{-1}$.
\end{proof}

From the above proposition, we observe that the bounds are tighter when the functions $(\phi_i')^{-1}$ have small Lipschitz constants. The closeness of the bounds also depends implicitly on $\epsilon$ via the Lagrange multipliers. 

Finally, we consider bounds for BW divergences with Mahalanobis distance generators.

\begin{proposition}[Mahalanobis distance] \label{prop:wc_signed_choquet_int_m_dist}
Let $I_h$ be a signed Choquet integral with representation \eqref{eq:quantile_rep} and distortion weight function $\gamma$ satisfying $0< \norm{\gamma}_2< \infty$. Let $g:\R^n\to\R$ be a $K$-Lipschitz function (wrt the $\L^2$ norm) with decomposition \eqref{eq:decomposition} and $\X\in\L^p_n$ with $g(\X)\in\L^p$. Consider $\bphi(\x)=\x^TQ\x$, where $Q$ is a positive-definite diagonal $n\times n$ matrix with smallest and largest eigenvalues $q_{\min}>0$ and $q_{\max}>0$, respectively. Then, for any $\epsilon>0$, 

\begin{enumerate}[label = $\roman*)$]
\item If $\gamma$ is non-decreasing, then 
\begin{equation*} I_h(g(\X))+\norm{\bbeta}_2\sqrt{\tfrac{\ep}{q_{\max}}}\norm{\gamma}_2\leq\sup_{\Y\in \bmfB^{n}_{\bphi, \epsilon}(\X)}I_h\big(g(\Y)\big)\leq I_h(g(\X))+K\sqrt{\tfrac{\ep}{q_{\min}}}\norm{\gamma}_2. 
\end{equation*}
\item If $\gamma$ is non-negative, then 
\begin{align*} 
\int_{0}^{1}\gamma(u)\iso{\Big(F_{g(x)}^{-1}(u)+\frac{1}{2\,q_{\max}\, \underline{\lambda}}\;\gamma(u)\Big)}du &\leq \sup_{\Y\in \bmfB^{n}_{\bphi, \epsilon}(\X)}I_h\big(g(\Y)\big) \\
&\leq \int_{0}^{1}\gamma(u)\iso{\Big(F_{g(x)}^{-1}(u)+\frac{1}{2\,q_{\min}\, \overbar{\lambda}}\;\gamma(u)\Big)}du, 
\end{align*}
where $\overbar{\lambda}$ is the unique positive solution to the equation in $\lambda$
\begin{equation*}
\int_{0}^{1}\Big(F_{g(\X)}^{-1}(u)-\iso{\big(F_{g(x)}^{-1}(u)+\frac{1}{2\,q_{\min}\,\lambda}\;\gamma(u)\big)}\Big)^2\,du=\tfrac{K^2\epsilon}{q_{\min}}\,,
\end{equation*}
and $\underline{\lambda}$ is the unique positive solution to the equation in $\lambda$
\begin{equation*}
\int_{0}^{1}\Big(F_{g(\X)}^{-1}(u)-\iso{\big(F_{g(x)}^{-1}(u)+\frac{1}{2\, q_{\max}\, \lambda}\;\gamma(u)\big)}\Big)^2\,du=\tfrac{\norm{\bbeta}_2^2\epsilon}{q_{\max}}\,.
\end{equation*}
\end{enumerate}
\end{proposition}

\begin{proof}
The proof is omitted as it follows similar steps as the proof of \Cref{prop:wc_signed_choquet_int}, with the difference that we apply \Cref{theorem:m_distance} instead of \Cref{prop:Lipschitz} in the first step.
\end{proof}

The bounds for signed Choquet integrals with uncertainty sets characterized by BW divergences with squared Mahalanobis distance generators take a similar form as those with the Wasserstein distance. The key difference is that for the Mahalanobis distance generator, the penalty term for the upper (resp. lower) bound is proportional to $K\sqrt{\frac{\epsilon}{q_{\min}}}$ (resp. $\norm{\bbeta}_2\sqrt{\tfrac{\epsilon}{q_{\max}}}$) rather than $K\epsilon$ (resp. $\norm{\bbeta}_2\epsilon$). The difference in the square root factor is because the Mahalanobis uncertainty set is defined using the squared Mahalanobis distance. We recover the bounds for the squared Wasserstein distance when $q_{\max}=q_{\min}=1$.

Since the $\L^2$ norm (\Cref{prop:wc_signed_choquet_int}) and the Mahalanobis distance (\Cref{prop:wc_signed_choquet_int_m_dist}) are separable Bregman generators, \Cref{prop:sep_Bregman} provides alternative upper and lower bounds to their worst-case risk when $\X$ is comonotonic. In \Cref{tab:diff_bounds}, we summarize the different bounds for the case when $\gamma$ is non-decreasing (i.e., $I_h$ is subadditive) and $g(\x)=\bbeta^T\x$ for some $\bbeta\in\R^n_{\geq}$. We note that \Cref{prop:wc_signed_choquet_int} and \Cref{prop:wc_signed_choquet_int_m_dist} yield tighter bounds, but \Cref{prop:sep_Bregman} applies to a larger class of BW divergences.

\begin{table}[htpb]
\centering
\caption{Different bounds for the Wasserstein and Mahalanobis distances.}\label{tab:diff_bounds}

\smallskip

\begin{tabular}{c c c}
Distance & Lower Bound & Upper Bound\\
\toprule
\toprule
Wasserstein (\Cref{prop:wc_signed_choquet_int}) &
$I_h(\bbeta^T\X)+\norm{\bbeta}_2\norm{\gamma}_2\epsilon$ &
$I_h(\bbeta^T\X)+K\norm{\gamma}_2\epsilon$\\
Wasserstein (\Cref{prop:sep_Bregman}) &
$I_h(\bbeta^T\X)+\tfrac{\norm{\bbeta}_1}{n}\norm{\gamma}_2\epsilon$ &
$I_h(\bbeta^T\X)+\norm{\bbeta}_1\norm{\gamma}_2\epsilon$\\
\midrule
\midrule
Mahalanobis (\Cref{prop:wc_signed_choquet_int_m_dist}) &
$I_h(\bbeta^T\X)+\norm{\bbeta}_2\sqrt{\tfrac{\epsilon}{q_{\max}}}\norm{\gamma}_2$ &
$I_h(\bbeta^T\X)+K\sqrt{\tfrac{\epsilon}{q_{\min}}}\norm{\gamma}_2$ \\
Mahalanobis (\Cref{prop:sep_Bregman}) &
\footnotesize{$I_h(\bbeta^T\X)+(\tfrac{1}{\sqrt{n}}\sum_{i\in\mN}\beta_iq_i^{-\frac{3}{2}})\sqrt{\epsilon}\norm{\gamma}_2$} &
\footnotesize{$I_h(\bbeta^T\X)+(\sum_{i\in\mN}\beta_iq_i^{-\frac{3}{2}})\sqrt{\epsilon}\norm{\gamma}_2$} \\
\bottomrule
\bottomrule
\end{tabular}
\end{table}

\section{Numerical Examples}\label{sec:examples}

In this section, we apply our results on bounding the DRO problem \eqref{eq:DRO-g-formulation} to an investment problem. We assume that the investor can choose to invest in $n=500$ different companies and for $i\in\mN$, the $i$-th company's stock has initial price $x^0_i$. Furthermore, for $i\in\mN$, the annual return on the $i$-th stock, denoted $R_i$, is given by 
\begin{equation*}
R_i=\mu_i+\sigma_i\,T_i,
\end{equation*}
where $\mu_i\in\R$, $\sigma_i>0$, and $T_i$ follows a $t$ distribution with 4 degrees of freedom. We further assume that the dependence structure of $(T_1, \ldots, T_{500})$ is given by  a $t$ copula with 4 degrees of freedom and correlation coefficient $0.7$. Thus, for $i\in\mN$, the mean annual return for the $i$-th stock is $\mu_i$, the variance is $2\,\sigma_i^2$, and the price in one year's time is $X_i:=x^0_i(1+R_i)$. Throughout, we assume that for all $i\in\mN$, $\mu_i$ and $\sigma_i$ are constant model parameters estimated by the investor and the current stock prices $\x^0:=(x^0_1, \ldots, x^0_{500})$ are known to the investor.

We further assume that the investor has initial capital $C$ and for each $i\in\mN$, the investor can choose to invest in $\beta_i$ units of the $i$-th stock or $\eta_i$ units of a European call or put option with strike price $K_i$ and maturity of 1 year. We denote by $c$ the difference between the initial capital and the initial value (cost) of the portfolio. Thus, $c>0$ represents uninvested cash and $c<0$ represents a loan. For simplicity, we assume that there is no interest. Therefore, the investor's loss in one year is given by
\begin{equation*}
g_m(\X)= C-\Big(\sum_{i=1}^{n-m}\beta_i\,X_i+\sum_{i=n-m+1}^{n}\eta_i\,\Psi_i(X_i, K_i)+c\Big),
\end{equation*}
where $0\leq m\leq n$ is the number of different options purchased by the investor and $\Psi_i(X_i, K_i)$ is the payoff of a European put or call option on stock $X_i$ with strike price $K_i$. Note that we assume without loss of generality that the investor chooses to invest in stock for the first $n-m$ companies. 

The investor aims to calculate lower and upper bounds to 
\begin{equation}\label{eq:DRO-example}
\sup_{\Y\in \bU^{n}_{\epsilon}(\X)}I_h(g_m(\Y))\,,    
\end{equation}
for different signed Choquet integrals, where $\bU^{n}_{\epsilon}(\X)$ is defined using the $p$-Wasserstein distance with norm $a$. By \eqref{eq:wass_bound_p}, there are four possible cases for the bounds, which we summarize in \Cref{tab:bounds}. Throughout, we let $\bzeta:=(\bbeta, \boldeta)$ and recall that $\tfrac{1}{p}+\tfrac{1}{q}=\tfrac{1}{a}+\tfrac{1}{b}=1$. 
\begin{table}[tpb]
\centering
\caption{Bounds for different choices of $a$ and $m$.}\label{tab:bounds}

\smallskip

\begin{tabularx}{\textwidth}{l l l L L}
$a$ & $\bbeta$ & $m$ & Lower Bound & Upper Bound\\
\toprule 
\toprule 
$a=1 \quad$ & $\norm{\bbeta}_{\infty}\geq\norm{\boldeta}_1 \quad$ & $0\leq m\leq n \quad$ & $I_h(g_m(\X))+\norm{\bbeta}_{\infty}\norm{\gamma}_q\epsilon$ & $I_h(g_m(\X))+\norm{\bbeta}_{\infty}\norm{\gamma}_q\epsilon$ \\

$a=1$ & $\norm{\bbeta}_{\infty}<\norm{\boldeta}_1$ & $1\leq m\leq n$ & $I_h(g_m(\X))+\norm{\bbeta}_{\infty}\norm{\gamma}_q\epsilon$ & $I_h(g_m(\X))+\norm{\eta}_{1}\norm{\gamma}_q\epsilon$ \\

$a\geq 2$ & $\bbeta\in\R^{n-m}$ & $m=0$ & $I_h(g_m(\X))+\norm{\bbeta}_{b}\norm{\gamma}_q\epsilon$ & $I_h(g_m(\X))+\norm{\bbeta}_{b}\norm{\gamma}_q\epsilon$\\

$a\geq 2$ & $\bbeta\in\R^{n-m}$ & $1\leq m\leq n$ & $I_h(g_m(\X))+\norm{\bbeta}_{b}\norm{\gamma}_q\epsilon$ & $I_h(g_m(\X))+\norm{\bzeta}_{b}\norm{\gamma}_q\epsilon$\\
\bottomrule
\bottomrule
\end{tabularx}
\end{table}
Note that the upper and lower bounds are equal in the first and third row of \Cref{tab:bounds}. Furthermore, when $m=0$, that is the aggregation function is linear, we recover proposition 2 $iv)$ from \citet{wozabal2014robustifying}. Additionally, when the upper and lower bounds do not coincide, the difference between the bounds is at most $\norm{\boldeta}_b\norm{\gamma}_q\epsilon$.

Next, we illustrate the quantile functions attaining the bounds in \Cref{tab:bounds} when $a=p=2$. For $m\in\{0, 1, 250, 500\}$, we assume that the investor invests $\$15$ in the first $n-m$ stocks (i.e., $\beta_i=\tfrac{15}{x_i^0}$, $i\in\{1, \ldots, n-m\}$) and purchases $3$ at-the-money call options for each of the $m$ remaining companies. To calculate the uninvested capital $c$, we assume the investor has initial capital $C=8000$, and the call options are priced using the Black-Scholes formula with risk free rate $0$ and implied volatility $0.2$. To generate our simulated stock data, we sample the initial stock value $x^0_i$ from a Uniform$(13,180)$ distribution, $\mu_i$ from a Uniform$(0.04,0.08)$ distribution, and $\sigma_i$ from a Uniform$(0.07,0.09)$ distribution. We further assume that stocks with a larger initial price will have larger mean returns and volatility. 

We first consider the ES and calculate the bounds for \eqref{eq:DRO-example}. Recall that the ES at level $\alpha\in(0,1)$ is $\text{ES}_{\alpha}(X)=\frac{1}{1-\alpha}\int_{\alpha}^{1}F_X^{-1}(t)\,dt$ and has distortion weight function $\gamma(u)=\frac{1}{1-\alpha}\Id_{u>\alpha}$. By \Cref{prop:wc_signed_choquet_int}, $i)$, noting that $\norm{\gamma}_2=(1-\alpha)^{-\frac{1}{2}}$ and $\gamma$ is non-decreasing, we have
    \begin{equation} \label{eq:es_bound}
\text{ES}_{\alpha}\big(g_m(\X)\big)+\norm{\bbeta}_2(1-\alpha)^{-\frac{1}{2}}\,\epsilon
    \leq
    \sup_{\Y\in \bU^{500}_{\epsilon}(\X)}\text{ES}_{\alpha}\big(g_m(\Y)\big)\leq\text{ES}_{\alpha}\big(g_m(\X)\big)+\norm{\bzeta}_2(1-\alpha)^{-\frac{1}{2}}\,\epsilon\,.
\end{equation}
Moreover, the quantile functions attaining the upper and lower bounds, respectively, are
\begin{subequations}
\begin{align} 
G_{u}^{-1}(u) &= F^{-1}_{g_{m}(\X)}(u)+\norm{\bbeta}_2(1-\alpha)^{-\frac{1}{2}}\,\epsilon\;\Id_{u>\alpha}, \label{eq:wc_es}\\
G_{l}^{-1}(u) &= F^{-1}_{g_{m}(\X)}(u)+\norm{\bzeta}_2(1-\alpha)^{-\frac{1}{2}}\,\epsilon\;\Id_{u>\alpha} \label{eq:wc_es2}.
\end{align}
\end{subequations}
Thus, the upper (resp. lower) bound in \eqref{eq:es_bound} is attained by applying a vertical shift of $\norm{\bbeta}_2(1-\alpha)^{-\frac{1}{2}}\,\epsilon$ (resp. $\norm{\bzeta}_2(1-\alpha)^{-\frac{1}{2}}\,\epsilon$) to the right tail of the reference quantile function. As $\alpha$ increases, the location of the shift moves to the right, and the magnitude of the shift (and hence the value of the worst-case ES) increases.

\begin{figure} [!htbp]
\centering
\includegraphics[width=0.4\textwidth]{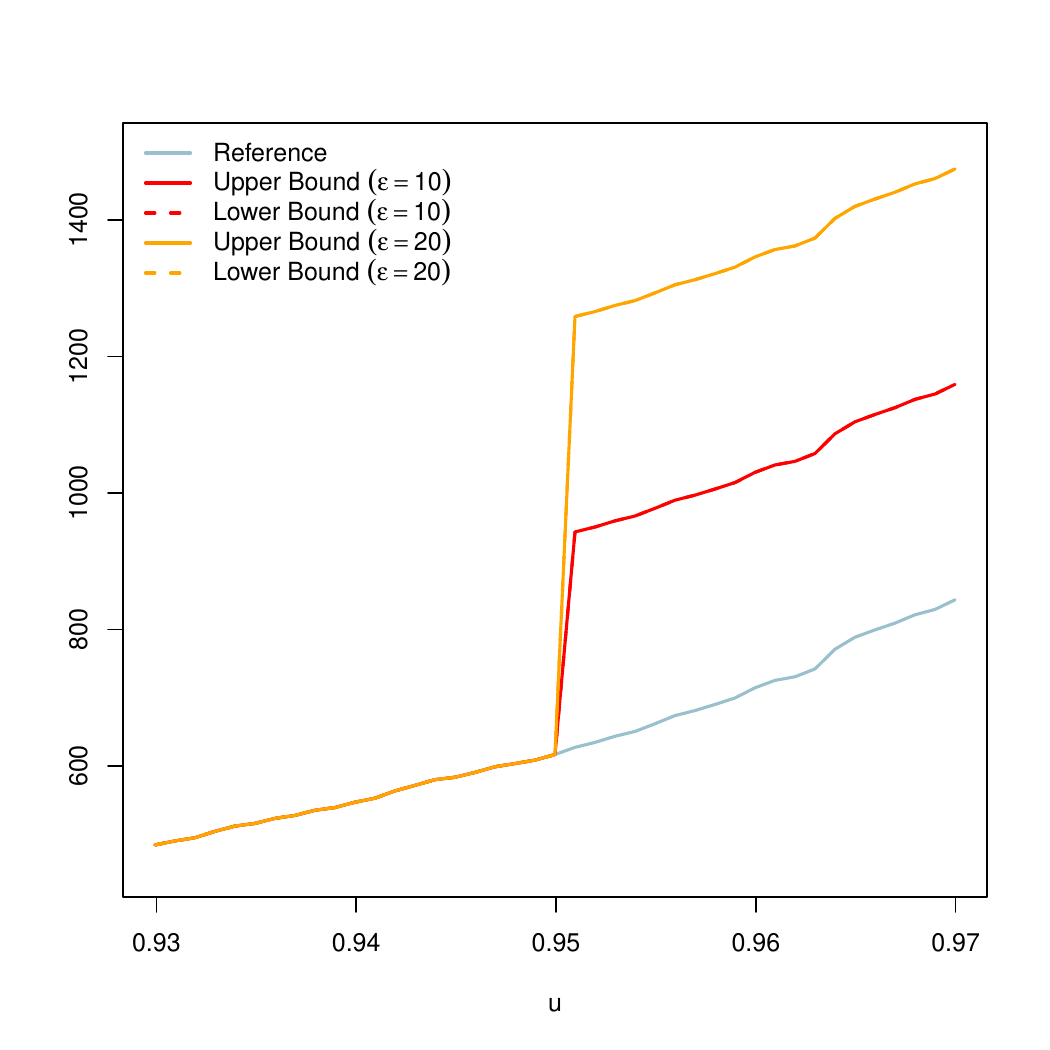}
\includegraphics[width=0.4\textwidth]{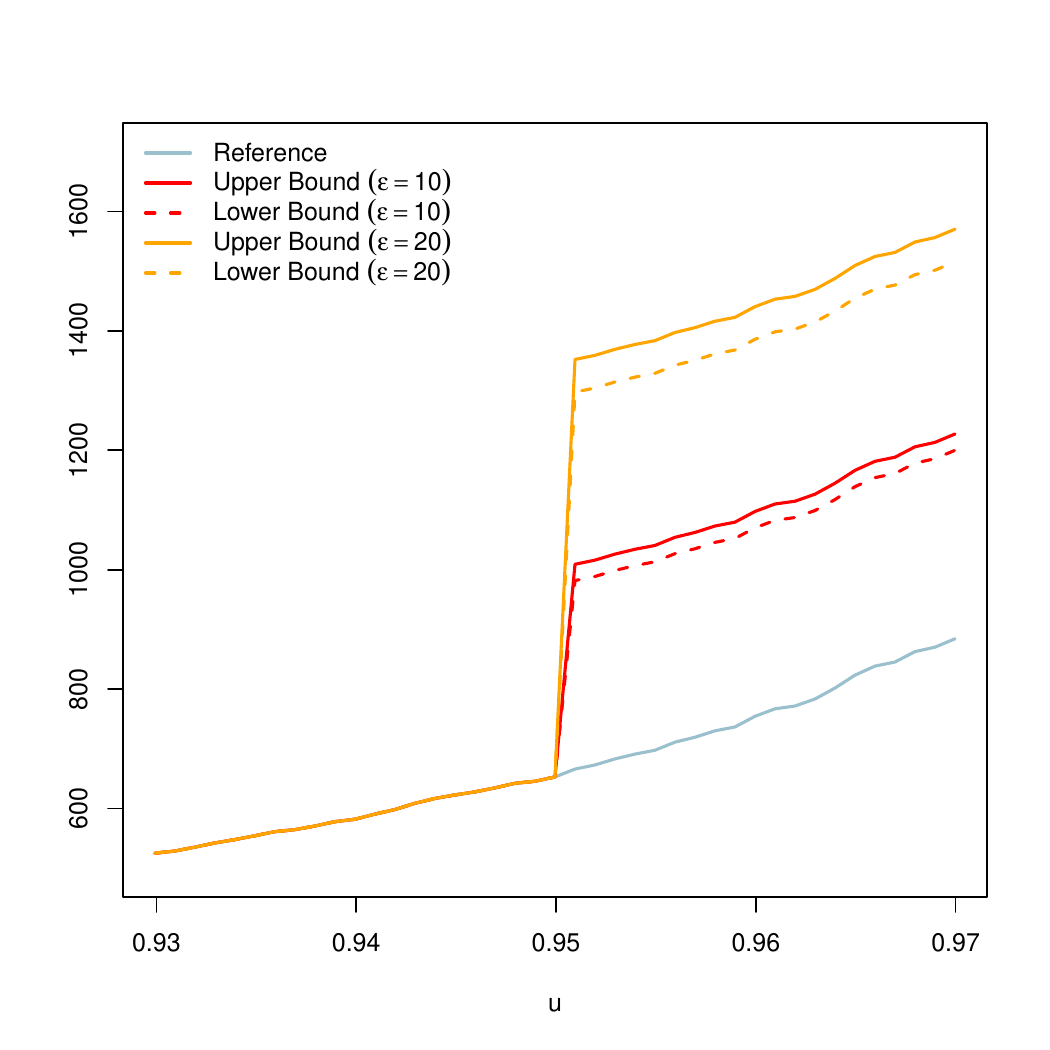}
\vspace{0.02\textwidth}
\includegraphics[width=0.4\textwidth]{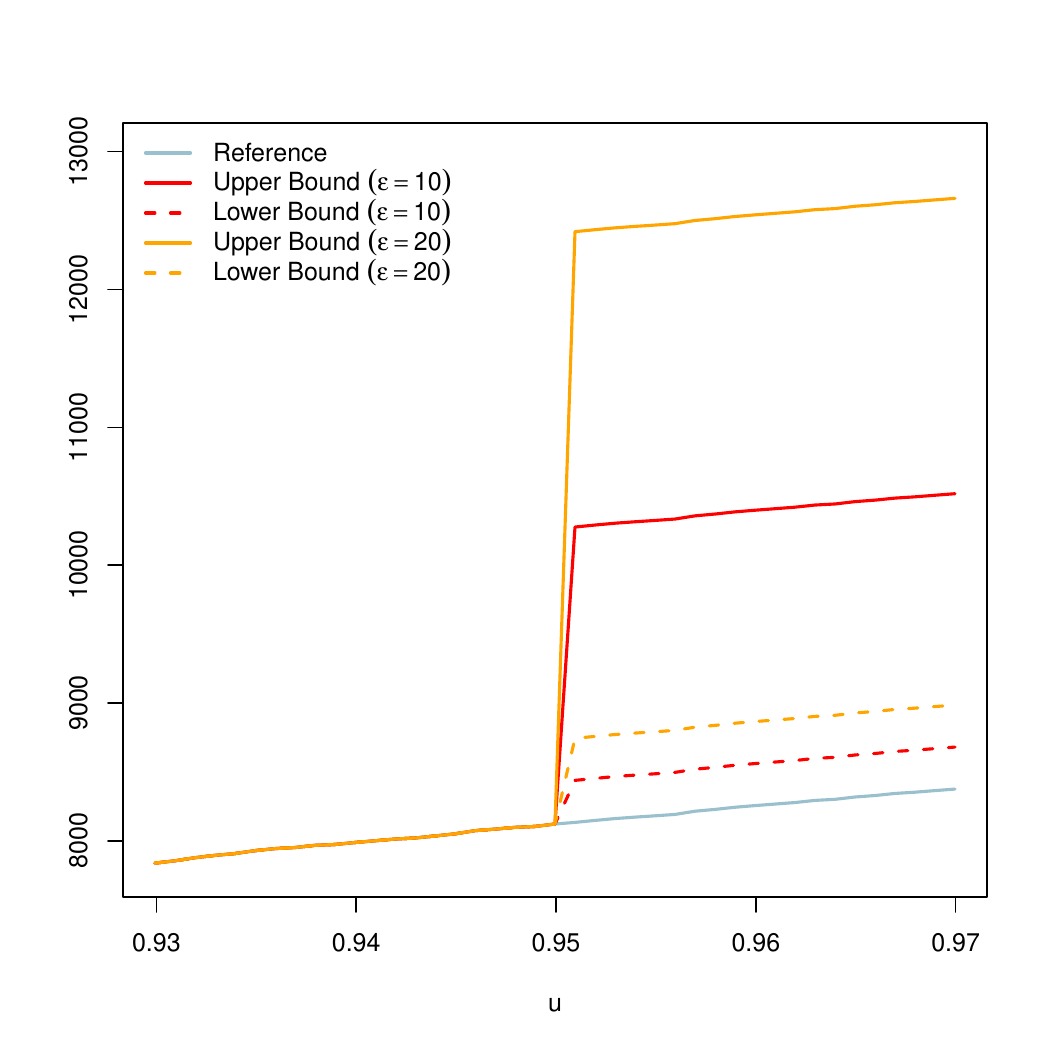}
\includegraphics[width=0.4\textwidth]{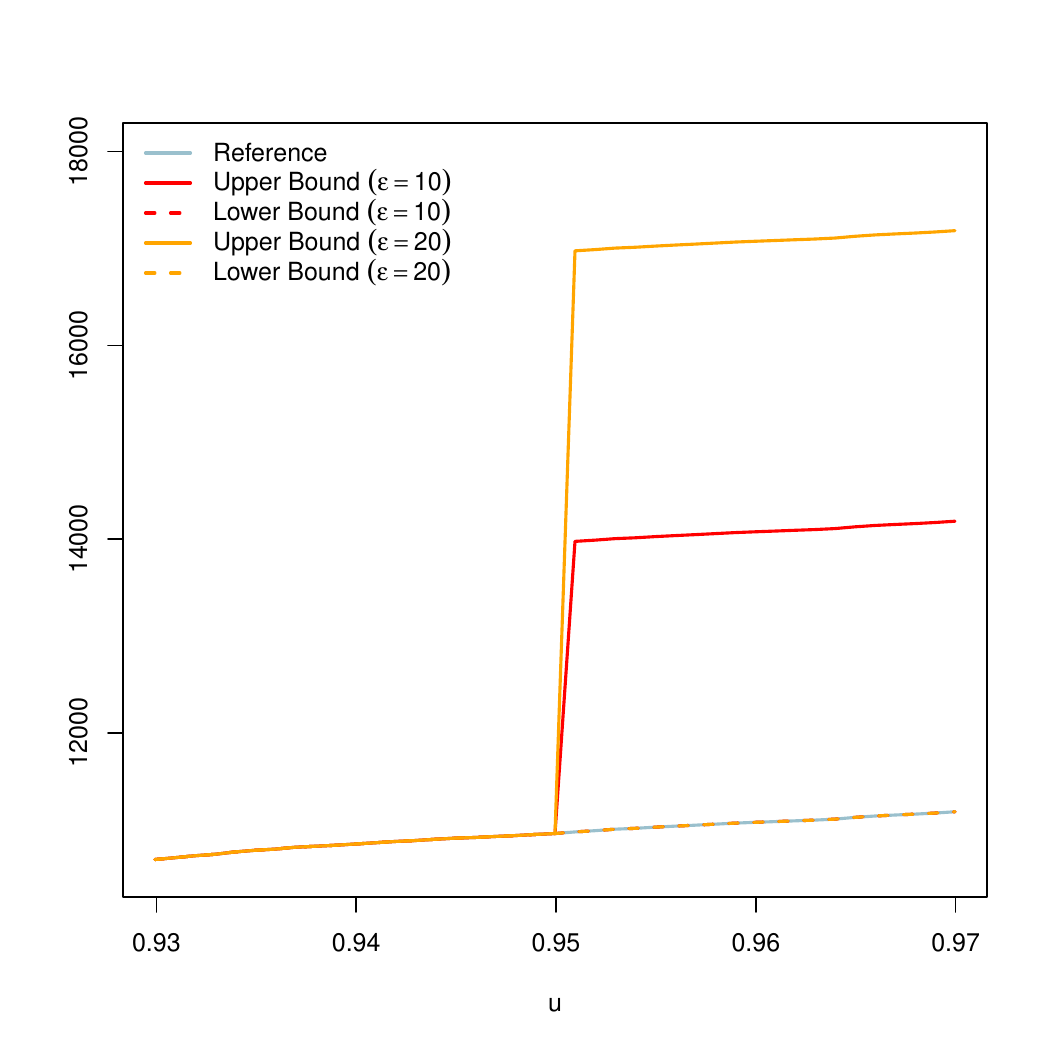}
\caption{Worst-case quantile functions for the ES (with $\alpha=0.95$). The reference distribution is plotted in blue. The red and orange lines correspond to the quantile functions attaining the upper and lower bounds for $\epsilon\in\{10, 20\}$, respectively. The quantile functions for $m\in\{0, 1, 250, 500\}$ are shown in the top left, top right, bottom left, and bottom right, respectively. In the top left the bounds coincide (linear aggregation) while in the bottom right the lower bound for the worst-case ES is the ES of the reference distribution (as $\norm{\bbeta} = 0$).}\label{fig:2}
\end{figure}

\Cref{fig:2} plots the functions \eqref{eq:wc_es} and \eqref{eq:wc_es2} for $\alpha=0.95$, $\epsilon\in\{10, 20\}$, and $m\in\{0, 1, 250, 500\}$. The blue curve shows the reference quantile function, which we estimate by sampling from the reference distribution $\X$. Furthermore, the red and orange curves correspond to $\epsilon=10$ and $\epsilon=20$, respectively, and we use solid curves for the quantile functions attaining the upper bounds and the dotted curves for the quantile functions attaining the lower bounds. Note that the plot only shows the right tails of the quantile functions as all curves are identical on $[0, 0.95)$. In the top left plot, corresponding to $m=0$, the aggregation function is linear, and hence the quantile functions attaining the upper and lower bounds coincide. In this case, we can exactly compute the worst-case ES, and the value of the worst-case ES agrees with proposition 2 $iv)$ from \citet{wozabal2014robustifying}. As $m$ increases from $1$ (top right) to $500$ (bottom right) and the aggregation function becomes less and less linear, the quantile functions attaining the bounds shift further apart and the bounds for the worst-case ES become looser. When $m=500$ (i.e., the aggregation function has no linear components), the lower bound on the worst-case ES equal to the ES of the reference distribution. Additionally, as we increase the tolerance distance $\ep$, the worst-case ES increases linearly in $\epsilon$.

Next, we consider the Inter-Expected Shortfall Range (IER) risk measure. Recall that the IER at level $\alpha\in (0.5,1)$ is given by 
\begin{equation*}
\text{IER}_{\alpha}(X):=\frac{1}{1-\alpha}\left(\int_{\alpha}^{1}F_X^{-1}(t) \, dt-\int_{0}^{1-\alpha}F_X^{-1}(t)\, dt\right).  
\end{equation*}
Note that the IER is a subadditive signed Choquet integral, but not a distortion risk measure. Indeed, its distortion weight function is $\gamma (u)=\frac{1}{1-\alpha}(\Id_{\alpha<u\leq 1}-\Id_{0\leq u\leq 1-\alpha})$ and thus is negative for $0\leq u\leq 1-\alpha$. As $\gamma$ is non-decreasing, by \Cref{prop:wc_signed_choquet_int}, $i)$ we have
\begin{equation*}
\text{IER}_{\alpha}\big(g_m(\X)\big)+\sqrt{\frac{2}{1-\alpha}}\norm{\bbeta}_2\,\epsilon\leq \sup_{\Y\in \bU^{500}_{\epsilon}(\X)}\text{IER}_{\alpha}\big(g_m(\Y)\big)\leq\text{IER}_{\alpha}\big(g_m(\X)\big)+\sqrt{\frac{2}{1-\alpha}}\norm{\bzeta}_2\,\epsilon, 
\end{equation*}
and the quantile functions attaining the upper and lower bounds, respectively, are
\begin{subequations}
    \begin{align} 
G_{u}^{-1}(u) &= F^{-1}_{g_{m}(\X)}(u)+\frac{\norm{\bzeta}_2}{\sqrt{2(1-\alpha)}}\;\epsilon\;(\Id_{\alpha<u\leq 1}-\Id_{0<u\leq 1-\alpha}), \label{eq:wc_ier}\\
G_{l}^{-1}(u) &= F^{-1}_{g_{m}(\X)}(u)+\frac{\norm{\bbeta}_2}{\sqrt{2(1-\alpha)}}\;\epsilon\;(\Id_{\alpha<u\leq 1}-\Id_{0<u\leq 1-\alpha}) \label{eq:wc_ier2}.
\end{align}
\end{subequations}

\begin{figure} [!htbp]
\centering
\includegraphics[width=0.4\textwidth]{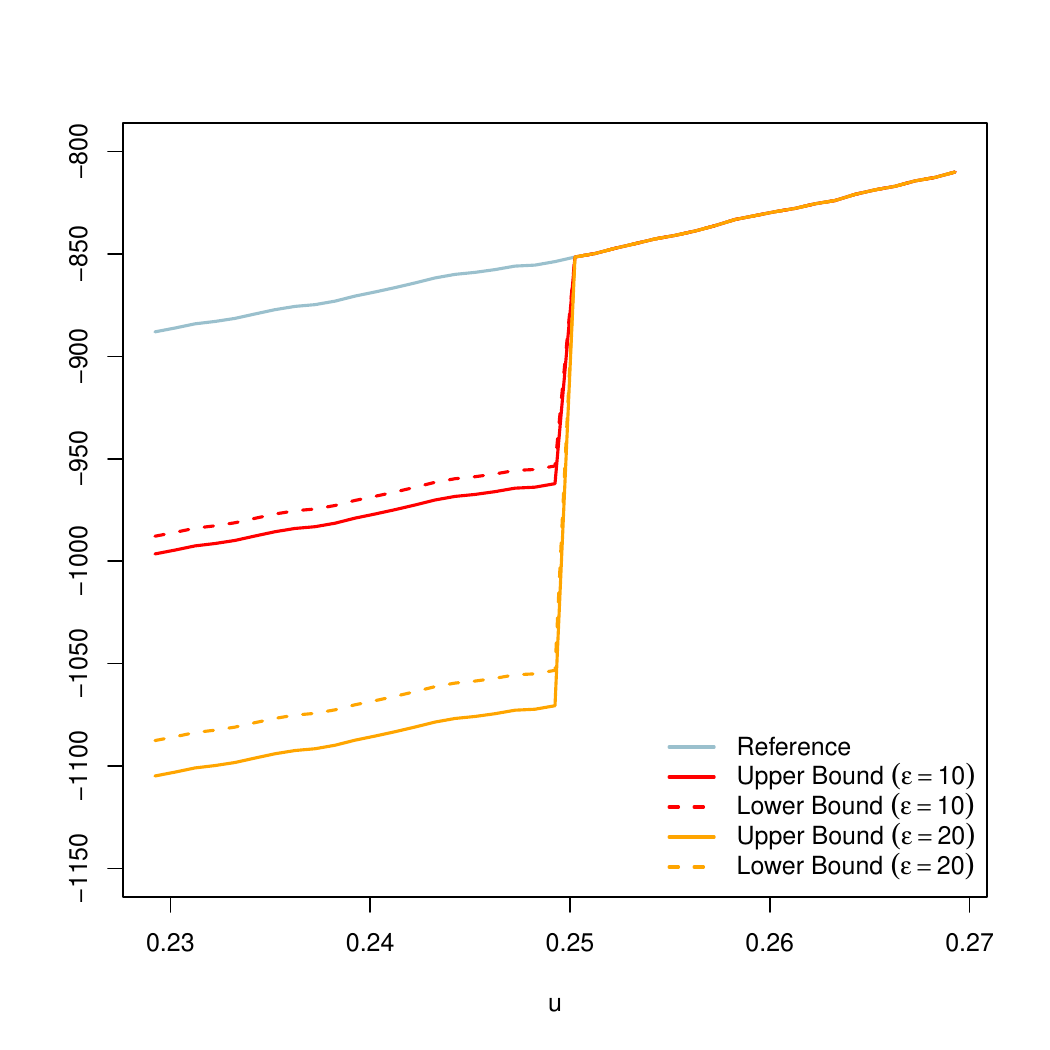}
\includegraphics[width=0.4\textwidth]{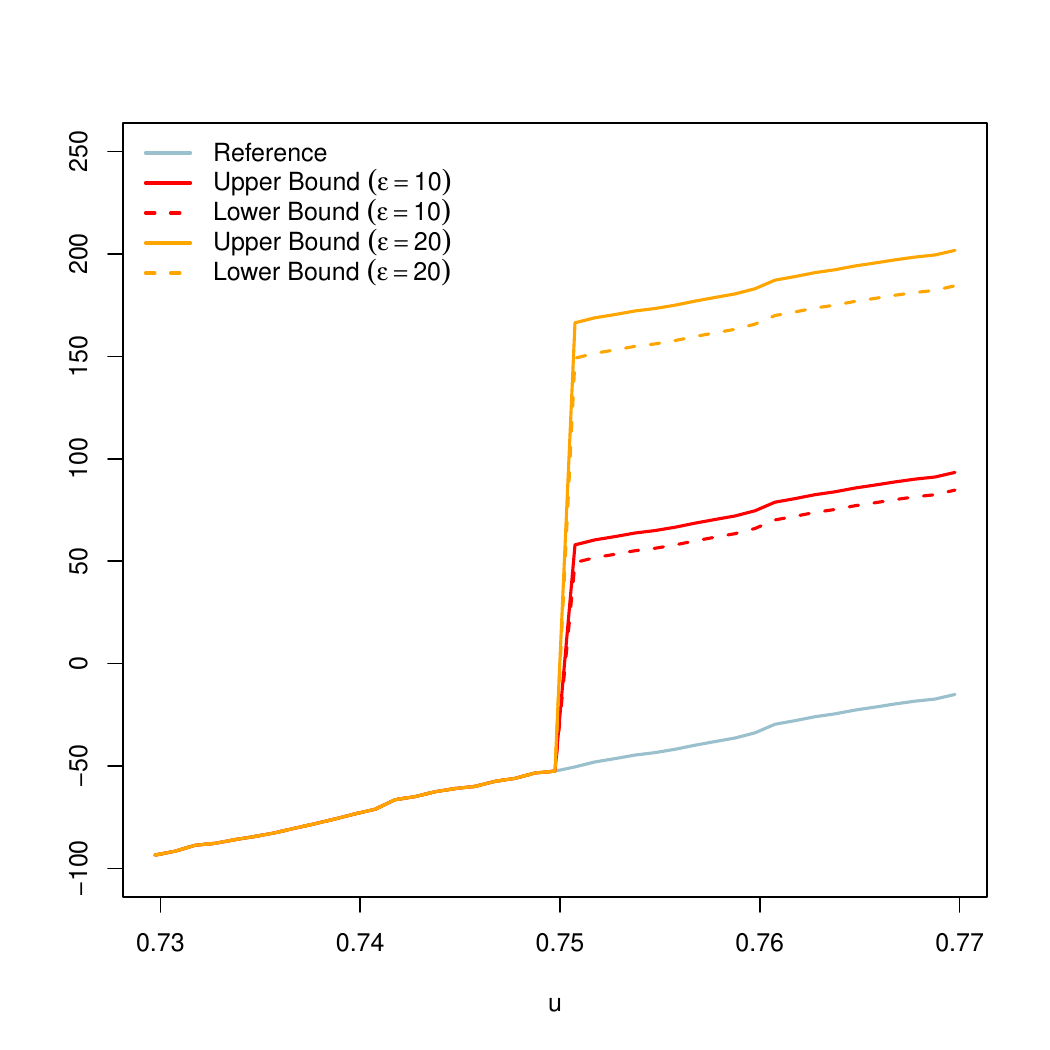}
\vspace{0.02\textwidth}
\includegraphics[width=0.4\textwidth]{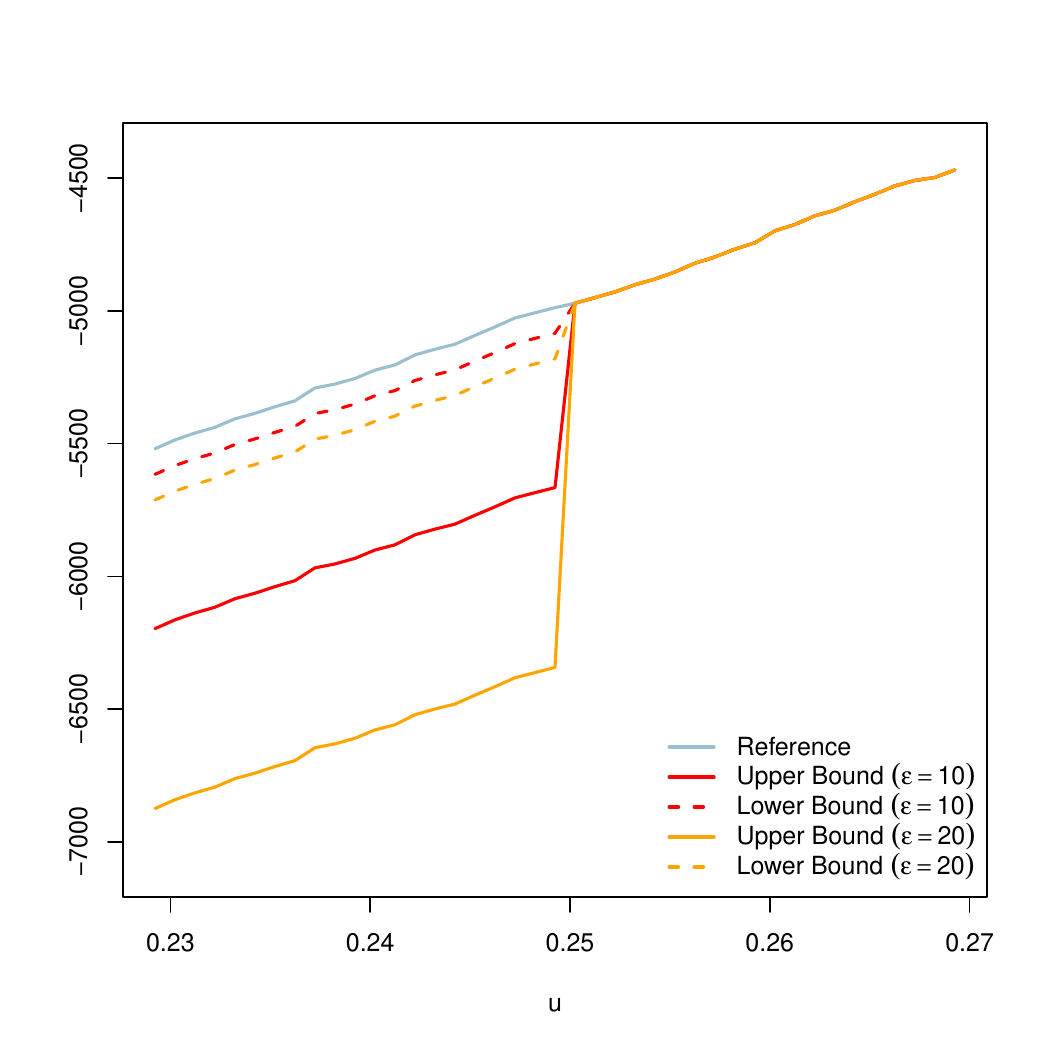}
\includegraphics[width=0.4\textwidth]{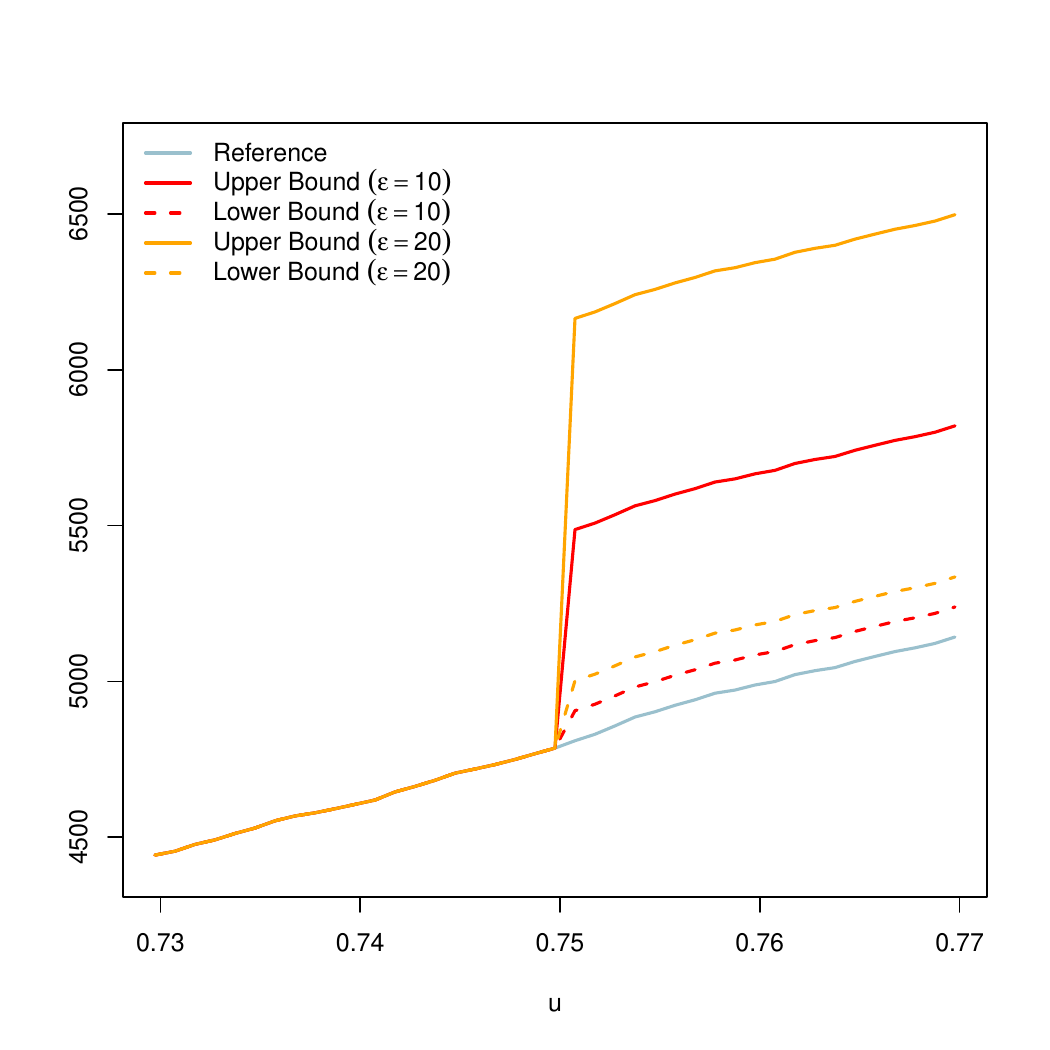}
\caption{Worst-case quantile functions for the IER (with $\alpha=0.75$). The reference distribution is plotted in blue. The red and orange lines correspond to the quantile functions attaining the upper and lower bounds for $\epsilon\in\{10, 20\}$, respectively. The quantile functions for $m\in\{1, 250\}$ are shown in the top and bottom plots, respectively.}\label{fig:3}
\end{figure}

\Cref{fig:3} displays the quantile functions \eqref{eq:wc_ier} and \eqref{eq:wc_ier2} for $\alpha=0.75$. We only plot $m = \{1, 250\}$, as similar to the ES, the upper and lower bounds coincide when the aggregation function is linear (i.e., $m=0$) and for $m = 500$, the lower bound for the worst-case IER is the IER of the reference distribution. The blue curve shows the reference quantile function, the red curves correspond to $\epsilon=10$, and the orange curves correspond to $\epsilon=20$. The solid red and orange curves correspond to the quantile functions attaining the upper bounds and the dotted curves correspond to the quantile functions attaining the lower bounds. We only display the quantile functions for $u\in [0.23, 0.27]$ (left) and $u\in [0.73, 0.77]$ (right) as all quantile functions are identical on $(0.25, 0.75)$. We observe that the quantile function of the lower bound is larger than the upper bound on $\{u\le0.25\}$ and smaller on $\{u > 0.75\}$. Similar to the ES example, the quantile functions attaining the bounds shift further apart as $m$ increases and the worst-case IER increases linearly in $\epsilon$. 

\begin{figure} [!htbp]
\centering
\includegraphics[width=0.4\textwidth]{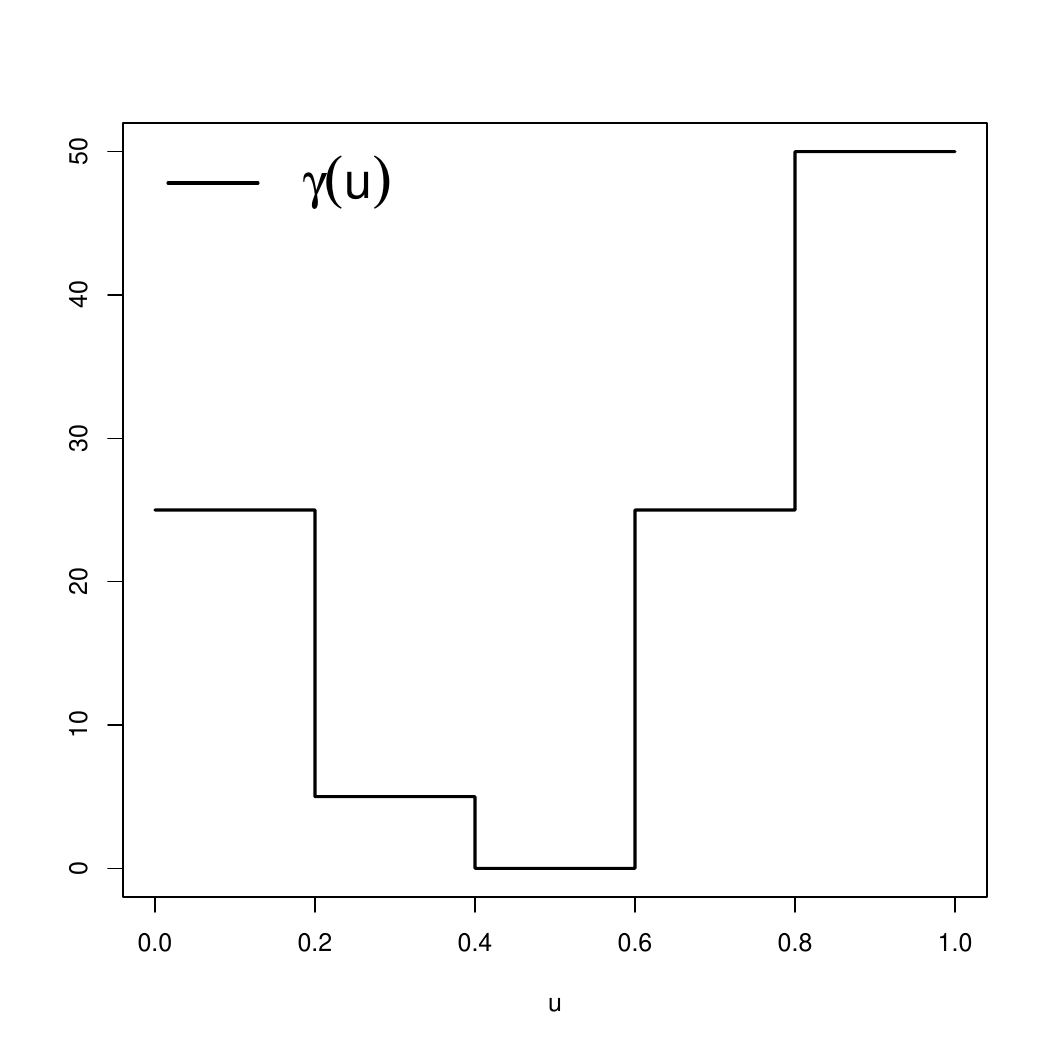}
\includegraphics[width=0.4\textwidth]{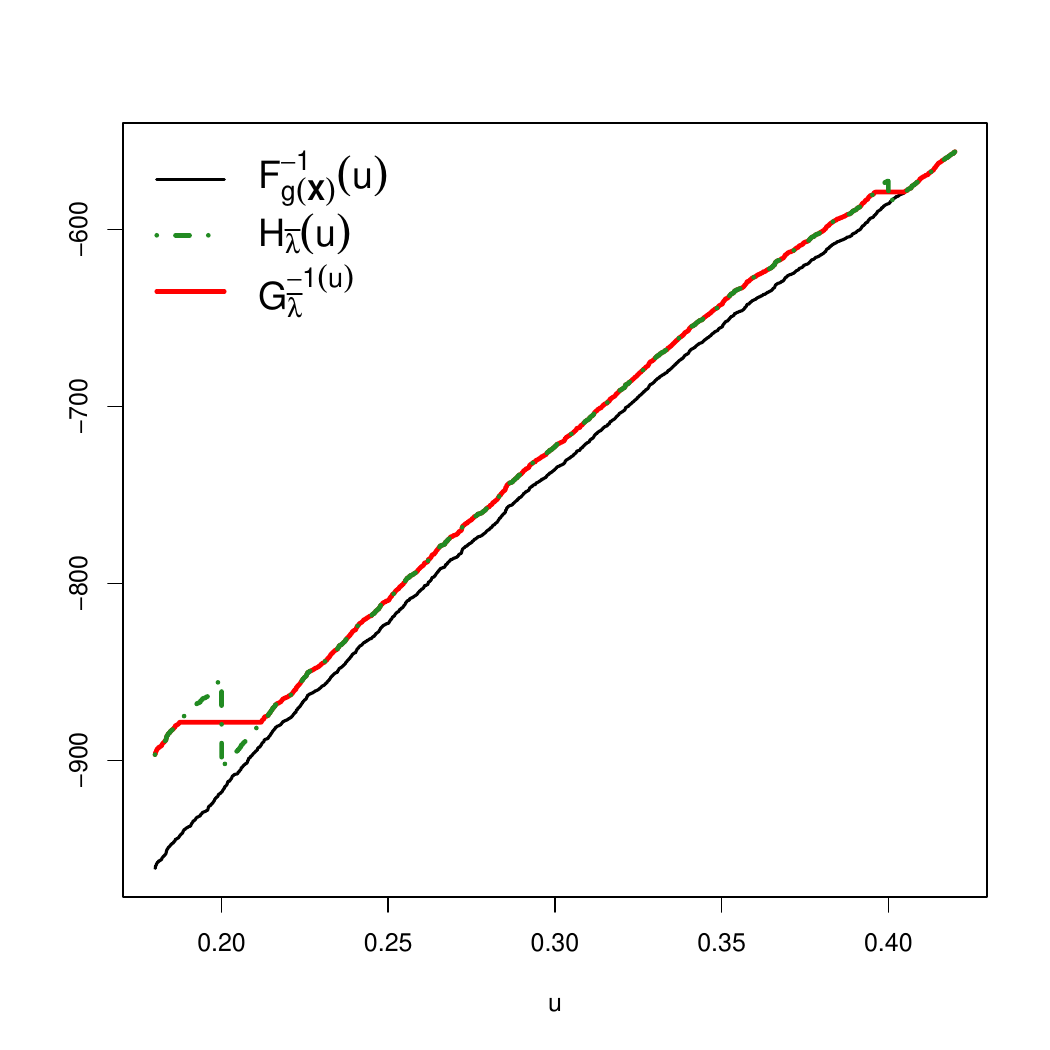}
\caption{Left: Inverse-S shaped distortion weight function $\gamma$. Right: Upper bound for worst-case quantile function $G_{\overbar{\lambda}}^{-1}$ (red, solid) for $\epsilon=10$ and $m=0$, reference quantile function $F^{-1}_{g_{0}(\X)}$ (black, solid), and the function $H_{\overbar{\lambda}}(u) = F^{-1}_{g_{0}(\X)}(u)+ \frac{\gamma(u)}{2\overbar{\lambda}} $ (green, dashed-dotted). Recall that $G_{\overbar{\lambda}}^{-1} =(H_{\overbar{\lambda}} )^\uparrow $.}\label{fig:4}
\end{figure}

Finally, we consider a piece-wise linear inverse-S shaped distortion weight function
\begin{equation*}
\gamma(u)=25\, (\Id_{0\leq u<0.2}+\Id_{0.6\leq u<0.8})+5(\Id_{0.2\leq u<0.4})+50(\Id_{0.8\leq u<1})\,,    
\end{equation*}
which is not a distortion risk measure as $\int_0^1 \gamma(u) \,du= 21 > 1$. We plot $\gamma$  in the left panel of \Cref{fig:4} and observe that the distortion weight function places the least weight at the centre of the distribution and the largest weight in the right tail. Since the aggregation function is the loss of a portfolio, this corresponds to weighting losses more than gains; a common practice in economics and often modelled by cumulative prospect theory (\citet{tversky}). 

To compute the bounds in \eqref{eq:wc_bound_2}, we require the isotonic projection of the function $H_{\lambda}(u):=F_{g_{m}(\X)}^{-1}(u)+\frac{1}{2\, \lambda}\, \gamma(u)$. We approximate the isotonic projection of $H_{\lambda}$ via the isotonic regression (as we evaluate $H_{\lambda}$ on a finite partition of $(0,1)$), which can be calculated using active set methods (\citet{leeuw}). We perform our computations in R using the function \texttt{gpava} from the package \texttt{isotone}.

The right panel of \Cref{fig:4} shows the worst-case quantile function with $\epsilon=10$ and $m=0$. The reference distribution in displayed in black. The green curve is $H_{\overbar{\lambda}}(u)$, and its isotonic projection is the worst-case quantile function $G^{-1}_{\overbar{\lambda}}$ plotted in red. We only display the functions for $u\in [0.18, 0.42]$ to better highlight the details in the graph. In particular, observe that the green curve is strictly decreasing around $u=0.2$ and $u=0.4$ and therefore is not a quantile function.

\begin{table}[!htpb]
\centering
\caption{Numerical values of the lower and upper bounds for different worst-case signed Choquet integrals and varying $m\in\{0,1,250, 500\}$. The last column corresponds to the difference between the upper and lower bound.}\label{tab:num_values}

\smallskip
\setlength{\tabcolsep}{10pt} 
\begin{tabular}{l r r r r r}
$I_h$ & $m$ & $I_h(g_m(\x))$ & lower bound & upper bound & length of bound\\
\toprule 
\toprule 
ES ($\alpha=0.95$) & 0 & 1,165.3 & 1,196.9 & 1,196.9 & 0\\
ES ($\alpha=0.95$) & 1 & 1,203.4 & 1,234.9 & 1,237.7 & 2.8\\
ES ($\alpha=0.95$) & 250 & 8,537.0 & 8,567.7 & 8,751.5 & 183.8\\
ES ($\alpha=0.95$) & 500 & 11,205.5 & 11,205.5 & 11,505.5 & 300.0 \\
\midrule
IER ($\alpha=0.75$) & 0 & 1,652 & 1,671.5 & 1,671.5 & 0\\
IER ($\alpha=0.75$) & 1 & 1,749 & 1,768.8 & 1,770.5 & 1.7\\
IER ($\alpha=0.75$) & 250 & 18,607 & 18,626.5 & 18,742.7 & 116.2 \\
IER ($\alpha=0.75$) & 500 & 23,840 & 23,839.9 & 24,029.6 & 189.7\\
\midrule
Piece-wise linear $\gamma$ & 0 & -3,578.2 & -3,384.3 & -3,384.3 & 0\\
Piece-wise linear $\gamma$ & 1 & -3,509.8 & -3,315.9 & -3,299.1 & 16.8\\
Piece-wise linear $\gamma$ & 250 & 23,157.4 & 23,344.7 & 24,474.0 &1,129.3\\
Piece-wise linear $\gamma$ & 500 & 45,568.7 & 45,568.7 & 47,411.8  & 1,843.1\\
\bottomrule
\bottomrule
\end{tabular}
\end{table}

Finally, we illustrate the numerical values of the bounds in \Cref{tab:num_values} for all three risk measures considered in this section for fixed $\epsilon=1$. As the bounds for ES and IER scale linearly with $\epsilon$, the value of the bounds for any $\epsilon\geq 0$ can be computed from the values in the table. Furthermore, we observe that the distance between the bounds increases as $m$ increases, and the bounds are best when the aggregation function is close to linear. Additionally, the distance between the bounds is very small relative to the absolute value of the reference risk. Indeed, for the worst bounds, that is for $m = 500$, the length of the bounds are less than 5\% of the reference risk measure, indicating that our method provides bounds for worst-case risks with small percentage error provided that $\epsilon$ is not too large.   

\section{Conclusion}

We prove that the image of multivariate Wasserstein uncertainty sets under Lipschitz aggregation functions are contained in and contain Wasserstein uncertainty sets of univariate rvs. This provides lower and upper bounds to risk-aware DRO problems with law-invariant risk functionals. We further show that for special cases of the aggregation function, the upper and lower bounds are equal, and therefore we can solve the multivariate DRO problem without directly computing the multivariate Wasserstein distance. Finally, we generalize our results to uncertainty sets characterized by asymmetric BW divergences and illustrate the bounds on the class of signed Choquet integrals.

\smallskip

\begin{APPENDICES}

\section{Isotonic Projections} \phantomsection\label{appendix:b}

Here, we prove some properties of the isotonic projection and show the max-min formulation of  isotonic projections, which is essential to prove \Cref{prop:iso_ordering}. The max-min formulation and its proof for the isotonic regression case can be found in \citet{barlow2}. 

\begin{proposition}[Properties] \label{prop:iso_properties}
Let $\ell\in\L^2((0,1)), \ k\geq 0$ and $c\in\R$. Then, all of the following hold:
\begin{enumerate} [label = $\roman*)$]
\item $\iso{(kl)}=k\iso{(\ell)}$,
\item $\iso{(\ell+c)}=\iso{\ell}+c$,
\item $\iso{\ell}(s)=\ell(s)+\sum_{i\in \I}(\theta_i-\ell(s))\Id_{s\in I_i}$,
where $\I$ is a countable index set, $I_i$ are mutually disjoint sub-intervals of $(0,1)$ with endpoints $a_i,b_i\in\R$ and $\theta_i\in\R$. Moreover, for any $i\in\I$,
\begin{equation} \label{eq:theta} 
\theta_i=\frac{1}{\abs{I_i}}\int_{a_i}^{b_i}\ell(s)ds,
\end{equation} 
where $\abs{I_i}:= b_i-a_i$ is the length of the interval $I_i$.
\end{enumerate}
\end{proposition}

\begin{proof}
Case $i)$ follows by Theorem 2.6 in \citet{brunk} since $\K$, the set of left-continuous and non-decreasing functions, is a closed convex cone.

Case $ii)$ follows by noting that
\begin{align*}
\iso{(\ell+c)} &= \argmin_{j\in\K}\int_{0}^{1}\big(j(s)-(\ell(s)+c)\big)^2ds 
=
\argmin_{j\in\K}\int_{0}^{1}\big((j(s)-c)-\ell(s)\big)^2ds\,.
\end{align*}
By definition of $\iso{\ell}$, the argmin of the right hand side $j^*$ satisfies $j^*-c=\iso{\ell}$, thus $\iso{(\ell+c)}=\iso{\ell}+c$. 

Case $iii)$, for a function $f$ we denote by $f^{*}$ its concave envelope. Then, as $\iso{\ell}(s)=\frac{d}{ds}\big(\int_{0}^{s}\ell(u)du\big)^{*}$, it holds that $\iso{\ell}(s)=\ell(s)+\sum_{i\in \I}(\theta_i-\ell(s))\Id_{s\in I_i}$ by Lemma 5.1 in \citet{brighi}. To prove \Cref{eq:theta}, recall that the intervals $I_i$ are disjoint, thus we have $\theta_i=\argmin_{\theta\in\R}\int_{I_i}(\theta-\ell(s))^2ds$. Calculating the argmin yields \eqref{eq:theta}.
\end{proof}

Before stating the max-min formulation for isotonic projections, we require additional definitions. 

\begin{definition} [Average Value of a Function] \label{definition:average} 
Let $S$ be a non-empty sub-interval of $(0,1)$ and $\ell\colon (0,1) \to \R$ be an integrable function. Then the average value of $\ell$ on $S$, is given by $\Av_l(S):= \frac{1}{\abs{S}}\int_Sl(s)ds$.
\end{definition}

For a set with only one element, i.e. $S=\{x\}$, it holds that $\Av_{\ell}(S)=\ell(x)$. Additionally, we drop the subscript $\ell$ of $\Av_{\ell}(\cdot)$, whenever it is clear from the context. 

\begin{definition} [Upper and Lower Sets - \citet{barlow2}, Def. 1.4.1] \label{definition:ul_sets} 
We define the following two notions of sets.
\begin{enumerate}[label = $\roman*)$]
    \item A set $L\subseteq (0,1)$ is a lower set wrt the quasi-order $\leq$, if for any $y\in L$ and any $x\in (0,1)$, the inequality $x\leq y$ implies that $x\in L$. 

    \item A set $U\subseteq (0,1)$ is an upper set wrt the quasi-order $\leq$, if for any $y\in U$ and any $x\in (0,1)$, the inequality $x\geq y$ implies that $x\in U$. 
\end{enumerate}
We denote the class of all lower sets by $S_{\L}$ and the class of all upper sets by $S_{\U}$.
\end{definition}

From the definition of a lower set, $S_{\L}$ is the set of all sub-intervals of $(0,1)$ with left endpoint $0$. Similarly, $S_{\U}$ is the set of all sub-intervals on $(0,1)$ with right endpoint $1$. Moreover, for any $\ell\in \K$ and $a\in\R$, the set $\{\ell\geq a\}:= \{s\in (0,1):\ell(s)\geq a\}$ is an upper set and $\{\ell\leq a\}:=\{s\in (0,1):\ell(s)\leq a\}$ is a lower set. This also holds when the inequalities in the sets are replaced with strict inequalities. 

We use the following lemma to prove the max-min formulation for isotonic projections. The citation for each part of the lemma references the corresponding result in \citet{barlow2} for isotonic regression. Note that cases $ii)$ and $iii)$ are proven in proposition A.3. of \citet{bernard}, but we include a short proof for completeness. 

\begin{lemma} \label{lemma:iso} 
Let $\ell\in\L^2((0,1))$. Then, all of the following hold:
\begin{enumerate}[label = $\roman*)$]
    \item (Theorem 1.3.6) For any function $\psi:\R\to\R$, $\int_{0}^{1}(\ell(s)-\iso{\ell}(s))\psi(\iso{\ell}(s))\, ds=0$. 
    \item (Equation 1.3.3) For any $f\in \K$, we have $\int_{0}^{1}(\ell(s)-\iso{\ell}(s))(\iso{\ell}(s)-f(s))\, ds\geq 0$.
    \item (Equation 1.3.8) For any $f\in \K$, we have $\int_{0}^{1}(\ell(s)-\iso{\ell}(s))f(s)\, ds\leq 0$.
    \item (Theorem 1.4.3) Let $a\in\R$, $L\in S_{\L}$ and $U\in S_{\U}$. Then, 
    \begin{enumerate} [label=$\alph*)$]
        \item $\Av(L\cap \{\iso{\ell}\geq a\})\geq a$,
        \item $\Av(L\cap \{\iso{\ell}> a\})> a$,
        \item $\Av(U\cap \{\iso{\ell}\leq a\})\leq a$,
        \item $\Av(U\cap \{\iso{\ell}< a\})< a$,
    \end{enumerate}
    whenever the sets are non-empty. 
\end{enumerate} 
\end{lemma}

\begin{proof}
\citet{barlow2} shows the above results for the isotonic regression, here we provide a proof for the isotonic projection. 

For case $i)$, let $\ell\in\L^2((0,1))$. Then, it holds that
\begin{align*}
\int_{0}^{1}\big(\ell(s)-\iso{\ell}(s)\big)\psi(\iso{\ell}(s))\,ds 
&= \sum_{i\in\I}\int_{I_i}(\ell(s)-\theta_i)\psi(\theta_i)\,ds 
= \sum_{i\in \I}\psi(\theta_i)\left(\int_{I_i}\ell(s)\, ds-|I_i|\theta_i\right) 
= 0,
\end{align*}
where the first and last equalities follow from \Cref{prop:iso_properties}, $iii)$. 

\smallskip

For case $ii)$, let $f\in \K$ and  $\alpha\in [0,1]$, then it holds that $(1-\alpha)\iso{\ell}+\alpha f \in \K$. 
By definition of the isotonic projection, the function $g(\alpha)=\int_{0}^{1}\big(\ell(s)-[(1-\alpha)\iso{\ell}(s)+\alpha f(s)]\big)^2\,ds$ attains its minimum on $[0,1]$ at $\alpha=0$. Moreover, $g'(0)=2\int_{0}^{1}\big(\ell(s)-\iso{\ell}(s)\big)\big(\iso{\ell}(s)-f(s)\big)\,ds$. Finally, since $g$ is quadratic in $\alpha$ and attains its minimum at $\alpha=0$, it follows that $g'(0)\geq 0$ and the desired result holds. 

\smallskip

For case $iii)$, we apply \Cref{lemma:iso}, $i)$ to the identity function and obtain $\int_{0}^{1}\big(\ell(s)-\iso{\ell}(s)\big)\iso{\ell}(s)\, ds=0$. Combining this fact with part $ii)$ of this lemma yields the inequality. 

\smallskip 

For case $iv)$, we only prove $b)$ as the remaining cases follow by similar arguments.
Let $a, a_1, a_2\in \R$ with $a_2>a_1$, and $\psi:\R\to\{0,1\}$ such that $\psi(x)=\Id_{a_1<x<a_2}$. By \Cref{lemma:iso}, $i)$, we have
\begin{equation} \label{eq:iso_4.1}
\int_{0}^{1}\big(\ell(s)-\iso{\ell}(s)\big)\psi(\iso{\ell}(s))\,ds=\int_0^1\big(\ell(s)-\iso{\ell}(s)\big)\Id_{a_1<\iso{\ell}(s)<a_2}\,ds=0\, .
\end{equation}
Since for any $U\in S_{\U}$, the function $\Id_{x\in U}$ is non-decreasing, we have \Cref{lemma:iso}, $iii)$, that
\begin{equation} \label{eq:iso_4.2}
\int_{0}^{1}\big(\ell(s)-\iso{\ell}(s)\big)\Id_{s\in U}\,ds\leq 0.
\end{equation} 
Therefore, for any $L\in S_{\L}$,
\begin{align}
\int_{L\cap\{\iso{\ell}>a\}}(\ell(s)-a)\,ds &> \int_{L\cap\{\iso{\ell}>a\}}\big(\ell(s)-\iso{\ell}(s)\big)\,ds 
\nonumber
\\
&= \int_{\{\iso{\ell}>a\}}\big(\ell(s)-\iso{\ell}(s)\big)\,ds - \int_{L^\complement\cap\{\iso{\ell}>a\}}\big(\ell(s)-\iso{\ell}(s)\big)\,ds\,,
\label{eq:iso_4-2-ineq}
\end{align}
where $L^\complement$ is the complement of $L$, and the equality follows since for any two sets $A$ and $B$, $A\cap B=B\setminus (A^\complement\cap B)$. Note that the first integral in \eqref{eq:iso_4-2-ineq} is zero by \Cref{eq:iso_4.1}. Furthermore, since $\iso{\ell}$ is non-decreasing, the set $L^\complement\cap\{\iso{\ell}>a\}$ is an upper set. Therefore, the second integral in \eqref{eq:iso_4-2-ineq} is non-positive by \Cref{eq:iso_4.2}. Combining these two facts yields, $\int_{L\cap\{\iso{\ell}>a\}}(\ell(s)-a)\,ds>0$. Rearranging the terms in the inequality gives 
\begin{equation*}
 a\,   < \, 
 \frac{1}{\abs{L\cap\{\iso{\ell}>a\}}}\int_{L\cap\{\iso{\ell}>a\}}\ell(s)\,ds   
 =
 \Av(L\cap\{\iso{\ell}>a\})\,,
\end{equation*}
provided $L\cap\{\iso{\ell}>a\}$ is non-empty. 
\end{proof}

We finally prove the max-min formulation of isotonic projections. 

\begin{lemma} (Max-Min Formulation)\label{lemma:max_min} Let $\ell\in\L^2((0,1))$ and fix $x\in (0,1)$. Then, 
\begin{enumerate}[label = $\roman*)$]
    \item $\displaystyle\max_{\{U\in S_{\U}|x\in U\}}\min_{\{L\in S_{\L}|x\in L\}}\Av(L\cap U)=\iso{\ell}(x)$, and
    \item $\displaystyle\max_{\{U\in S_{\U}|x\in U\}}\min_{\{L\in S_{\L}|L\cap U\neq\emptyset\}}\Av(L\cap U)=\iso{\ell}(x)$.
\end{enumerate} 
\end{lemma}

\begin{proof}
We only prove the second case, as the first case follows using similar arguments.

For $x\in (0,1)$, define $a:=\iso{\ell}(x)$ and $U_a:=\{\iso{\ell}\geq a\}$. For any $L \in S_{\L}$, we have that by \Cref{lemma:iso}, $iv)$ $a)$, $\Av(L\cap U_a\})\geq a$ whenever $L\cap U_a\neq \emptyset$. Therefore, 
\begin{equation*}
\min_{\{L\in S_{\L}|L\cap U_a\neq\emptyset\}}\Av(L\cap U_a)\geq a.
\end{equation*}
Additionally, as $x\in U_a$ it follows that 
\begin{equation} \label{eq:iso_2.1}
\max_{\{U\in S_{\U}|x\in U\}}\min_{\{L\in S_{\L}|L\cap U\neq\emptyset\}}\Av(L\cap U)\geq a.
\end{equation}
Next, let $U\in S_{\U}$ such that $x\in U$. Then, $\{\iso{\ell}\leq a\}\cap U$ is non-empty and $\Av(\{\iso{\ell}\leq a\}\cap U)\leq a$ by \Cref{lemma:iso}, $iv)$ $c)$. Therefore, 
\begin{equation} \label{eq:iso_2.2}
\max_{\{U\in S_{\U}|x\in U\}}\min_{\{L\in S_{\L}|L\cap U\neq\emptyset\}}\Av(L\cap U)\leq a
\end{equation}
Combining Equations \eqref{eq:iso_2.1} and \eqref{eq:iso_2.2} completes the proof. 
\end{proof}

\smallskip 

{\parindent0pt
\textbf{Proof of \Cref{prop:iso_ordering}:}
Let $\ell_1, \ell_2\in\L^2((0,1))$ satisfying $\ell_2(s)\leq \ell_1(s)$, for any non-empty interval $S\subseteq (0,1)$. Then, for any non-empty interval $S\subseteq (0,1)$ it holds 
\begin{equation} \label{eq:ordering_1}
\frac{1}{\abs{S}}\int_Sl_2(s)\,ds=\Av_{\ell_2}(S)
\leq
\Av_{\ell_1}(S)=\frac{1}{\abs{S}}\int_Sl_1(s)\,ds.
\end{equation}
Thus, by \Cref{lemma:max_min}, $i)$ we have for any fixed $s\in (0,1)$, that
\begin{equation} \label{eq:ordering_2}
\iso{\ell}_2(s)=\max_{\{U\in S_{\U}|x\in U\}}\min_{\{L\in S_{\L}|x\in L\}}\Av_{\ell_2}(L\cap U)\leq \max_{\{U\in S_{\U}|x\in U\}}\min_{\{L\in S_{\L}|x\in L\}}\Av_{\ell_1}(L\cap U)=\iso{\ell}_1(s).
\end{equation}

If $\ell_2(s)<\ell_1(s)$, then the inequalities in \eqref{eq:ordering_1} and \eqref{eq:ordering_2} hold with strict inequality. 
$\hfill\square$
}

\end{APPENDICES}

\bibliographystyle{informs2014}
\nocite{*}
\bibliography{Refs.bib}

\end{document}